\documentclass[draft]{amsart}
\usepackage{amsfonts,amssymb,latexsym,amsmath, amsxtra,url, mathrsfs}
\usepackage[utf8]{inputenc}

\usepackage{graphicx}
\graphicspath{ {./Figures/} }

\usepackage{enumerate}

\usepackage{tikz}
\usepackage{caption}

\usepackage{hyperref}
\hypersetup{colorlinks=true,linkcolor=blue }

\theoremstyle{plain}
\newtheorem{theorem}{Theorem}[section]
\newtheorem*{theorem*}{Theorem}
\newtheorem{corollary}[theorem]{Corollary}
\newtheorem{conjecture}[theorem]{Conjecture}
\newtheorem{lemma}[theorem]{Lemma}

\newtheorem*{conjecture*}{Conjecture}

\theoremstyle{definition}
\newtheorem{definition}[theorem]{Definition}

\theoremstyle{remark}

\numberwithin{equation}{section}

\newcommand{\Z}{\mathbb{Z}}

\renewcommand{\P}{\mathcal{P}}

\title[Mod 11 and 13 $A_2$ Rogers--Ramanujan Type Identities]{Proofs of Modulo 11 and 13 Cylindric Kanade-Russell Conjectures for $A_2$ Rogers--Ramanujan Type Identities}
\author[Ali Kemal Uncu]{Ali Kemal Uncu}
\address{Johann Radon Institute for Computational and Applied Mathematics, Austrian Academy of Science, Altenbergerstraße 69, A-4040 Linz, Austria}
\email{akuncu@ricam.oeaw.ac.at}
\address{University of Bath, Faculty of Science, Department of Computer Science, Bath, BA2 7AY, UK}
\email{aku21@bath.ac.uk}

\thanks{Research of the author is partly supported by EPSRC grant number EP/T015713/1 and partly by FWF grant P-34501N}

  

\date{\today}

\begin{document}

\maketitle

\begin{abstract} 
We present proofs of two new families of sum-product identities arising from the cylindric partitions paradigm. Most of the presented expressions, the related sum-product identities, and the ingredients for the proofs were first conjectured by Kanade--Russell in the spirit of Andrews--Schilling--Warnaar identities of the $A_2$ Rogers--Ramanujan type. We follow the footsteps of Kanade--Russell while we alter the computations heavily to accomplish our goals.
\end{abstract}

\section{Introduction}

There is an ever-growing synergy between number theory, combinatorics, $q$-series, and affine Lie algebras that led to groundbreaking techniques and beautiful mathematical discoveries. Among these are the Rogers--Ramanujan type identities where an infinite $q$-series is equal to a infinite product with a modular structure. First appeared at the intersection of number theory and combinatorics, the Rogers--Ramanujan identities have been of great interest. These sum-product identities have been studied, proved and generalized in many different ways over the years \cite{And74, And89, Bre79, Bre83, Cor17, GM81, Gor61, GOW16, Pas20}. These identities also naturally arose in many other fields including mathematical physics \cite{Bax81}, representation theory of affine Lie algebras and vector operator algebras \cite{LW84, LW85}, knot theory in relation to the colored Jones polynomials \cite{AD11}, and algebraic geometry \cite{BMS13} over the years. 

For some non-negative integer $L$ and formal variables $a$ and $q$, let $q$-Pochhammer symbol be $(a;q)_L := (1-a)(1-aq)\dots(1-aq^{L-1})$, and $(a;q)_\infty := \lim_{L\rightarrow\infty} (a;q)_L$, $\theta(a;q) := (a,q/a;q)_\infty$, and for $a_1,\dots,a_k$ some formal variables, define the shorthand notation $\theta(a_1,a_2,\dots, a_k;q):= \theta(a_1;q)\theta(a_2;q)\dots\theta(a_k;q)$. 

The Rogers--Ramanujan identities are as follows \cite{RR19}. \begin{theorem}[Rogers--Ramanujan identities]  \begin{equation}\label{eq:RR}\sum_{n\geq 0} \frac{q^{n^2}}{(q;q)_n} = \frac{1}{\theta(q;q^5)}\ \ \text{and}\ \ \sum_{n\geq 0} \frac{q^{n^2+n}}{(q;q)_n} = \frac{1}{\theta(q^2;q^5)}.\end{equation}
\end{theorem} The reciprocal $q$-Pochhammer products on the right-hand side of \eqref{eq:RR} has the $\pm 1$ and $\pm 2$ residue classes modulo 5, respectively. We call these \textit{modulo 5} identities.

A \textit{composition} $c$ of $n$ is a finite list of non-negative integers that sum up to $n$. A \textit{partition} is a composition where no element of the list (called \textit{parts}) are zero and the list elements are ordered in a non-increasing order. We define the size of a composition $\pi$ as the sum of all its parts and denote this by $|\pi|$. We denote the number of parts in a composition $\pi$ by $\#(\pi)$. A composition (resp. partition) with size $n$ is called ``\textit{a composition (resp. partition) of }$n$." The empty list is considered as the unique composition/partition of $0$ with 0 parts. For example, $(2,0,2)$ is a composition with 3 parts and $(1,1,1,1)$, $(4,3,1)$, and $(2,2)$ are partitions of 4, 8, and 4, respectively. 

MacMahon \cite{Mac16} and Schur \cite{Sch17} gave combinatorial interpretations to Rogers--Ramanujan identities independently.

\begin{theorem}[Combinatorial interpretaton of Rogers--Ramanujan identities]
\label{th:RRcomb}
Let $i=0$ or $1$. For every natural number $n$, the number of partitions of $n$ such that the difference between two consecutive parts is at least $2$ and the the smallest part is strictly greater than $i$ is equal to the number of partitions of $n$ into parts congruent to $\pm (1+i) \mod 5.$
\end{theorem}

Gordon \cite{Gor61} presented a wide generalization of Theorem \ref{th:RRcomb} to all odd modulus $\geq 5$.
\begin{theorem}[Gordon's identities, 1961]
\label{th:Gordon}
Let $r$ and $i$ be integers such that $r\geq 2$ and $1\leq i \leq r.$ The number of partitions $\pi=(\pi_1,\pi_2,\dots,\pi_s)$ of $n$ such that $\pi_{j}-\pi_{j+r-1} \geq 2$ for all $j$ with at most $i-1$ 1s appears as parts in $\pi$ are equal to the number of partitions of $n$ whose parts are not congruent to $0,\pm i \mod  2r+1$.
\end{theorem}
The Rogers--Ramanujan identities correspond to the cases $r=i=2$ and $r=2$, $i=1$.

Andrews found the $q$-series counterpart to Gordon's identities \cite{And74}.
\begin{theorem}[Andrews--Gordon identities, 1974]
\label{th:AGseries}
Let $r \geq 2$ and $1 \leq i \leq r$ be two integers. We have
\begin{equation}\label{eq:AGri}
\sum_{n_1\geq\dots\geq n_{r-1}\geq0}\frac{q^{n_1^2+\dots+n_{r-1}^2+n_{i}+\dots+n_{r-1}}}{(q;q)_{n_1}} {n_1\brack n_1-n_2}_q \cdots {n_{r-2}\brack n_{r-2}-n_{r-1}}_q=\frac{\theta(q^{i};q^{2r+1})(q^{2r+1};q^{2r+1})_\infty}{(q;q)_\infty},
\end{equation}
where for two integers $n$ and $m$,
$${m+n \brack m}_q := \begin{cases} 
\displaystyle
\frac{(q;q)_{m+n}}{(q;q)_m(q;q)_{n}}&\text{for }m, n \geq 0,
\vspace{3pt}
\\  \ 0&\text{otherwise,}\end{cases}$$
is the classical $q$-binomial coefficient.
\end{theorem}
Note that the Rogers--Ramanujan identities are the particular case of~\eqref{eq:AGri} where $r=i=2$, and $r=2$ and $i=1$. Interested readers can get a great overview of the history of the Rogers--Ramanujan identities, their significance, and some generalizations in the recent book of Sills \cite{Sil17}. 

The identities \eqref{eq:AGri} can be proven by the Bailey machinery coming from the world of $q$-series. This powerful mechanism starts with a pair of $q$-expressions, called a \textit{Bailey pair}, that satisfies a pre-defined relation and modifies this pair iteratively (using \textit{Bailey lemma} or one of its generalizations) to make a new Bailey pair (see \cite{AAB87,And86,Bai49, Sil17}). That way, by starting with the pair related to Rogers--Ramanujan identities, a whole infinite chain of identities  \eqref{eq:AGri} can be acquired.  The identities \eqref{eq:AGri} are certain characters related to affine Lie algebra $A_1^{(1)}$, and we thus refer to them as $A_1$ Rogers--Ramanujan identities. The original Bailey mechanism was later extended to $A_{n-1}$ for general $n$ \cite{ML92, ML95}. However,these works did not yield $A_{n-1}$ Rogers--Ramanujan identities. 

In their influential paper, Andrews, Schilling and Warnaar \cite{ASW99} were able to describe an $A_2$ Bailey lemma and the associated Bailey machinery. They found several infinite families of identities, One of their modulo 7 identities is as follows. 
\begin{theorem}[Andrews--Schilling--Warnaar, 1999]
\label{th:ASW}
\begin{equation}\label{eq:ASW400}\sum_{r_1,s_1\geq 0} \frac{q^{r_1^2-r_1 s_1+s_1^2+r_1+s_1}}{(q;q)_{r_1}} {2 r_1\brack s_1}_q=\frac{1}{\theta(q^2,q^3,q^3;q^7)}.\end{equation}
\end{theorem}

Andrews--Schilling--Warnaar found several very general families of sum-product identities. Of particular interest to representation theory, the product-sides of these identities are character formulas of the $W_3$ algebra multiplied by an extra factor $(q;q)^{-1}_\infty$ \cite{FW16}. These formulas do not yield manifestly positive sum-sides for the character formulas because of this extra factor. 

For example, one of Andrews--Schilling--Warnaar's modulo 10 identities after clearing the extra factor $(q;q)^{-1}_\infty$ is as follows.
\begin{theorem}[Andrews--Schilling--Warnaar, 1999]
\label{th:ASW500}
\begin{equation}\label{eq:ASW500}(q,q)_{\infty} \sum_{\substack{r_1\geq r2\geq 0\\s_1\geq s_2\geq 0}} \frac{q^{r_1^2- r_1 s_1  + s_1^2 + r_2^2 - r_2 s_2 + s_2^2   + r_1 + r_2 + s_1 + s_2}}{(q;q)_{r_1-r_2}(q;q)_{r_2}(q;q)_{s_1-r_2}(q;q)_{s_2}(q;q)_{r_2+s_2+1}} =\frac{1}{\theta(q^2,q^3,q^3,q^4,q^4,q^5;q^{10})}\end{equation}
\end{theorem}
Recall the Euler's Pentagonal Number Theorem \cite{And84b} \begin{equation}\label{EPNT}(q,q)_\infty = \sum_{i=-\infty}^\infty (-1)^i q^{i(3i+1)/2}.\end{equation}
Although it is easy to see that the right-hand side of \eqref{eq:ASW500} has positive coefficients, in light of \eqref{EPNT} this is not directly visible on the left-hand side. In contrast, both sides of \eqref{eq:ASW400} are manifestly positive. The manifestly positive sum representations give insight to the structure of certain modules for the affine Lie algebra $A_2^{(1)}$. These mentioned character of standard modules for the affine Lie algebra $A_2^{(1)}$. Interested readers can find more on this connection in \cite{ASW99, KR, KanF, LW84, LW85}.

Recently, the discovery of manifestly positive identities of these character formulas through a scheme with combinatorial roots attracted the attention and led to many new Rogers--Ramanujan type identities.

In 1997, Gessel and Krattenthaler \cite{GK97} defined \textit{cylindric partitions} in context of non-intersecting lattice paths. Borodin \cite{Bor07} gave univariate product formulas for the generating functions of the number of cylindric partitions. Foda and Welsh \cite{FW16} proved the $A_1$ Rogers--Ramanujan identities using the combinatorics of cylindric partitions. This led to Corteel's combinatorial proof of the Rogers--Ramanujan identities using cylindric partitions \cite{Cor17}. In 2019, Corteel and Welsh \cite{CW19} derived functional equations for the bivariate generating functions for the number the number of cylindric partitions using the largest part statistic. While doing so, they also gave a new proof of Andrews--Schilling--Warnaar's modulo 7 $A_2$ Rogers--Ramanujan identities (including \eqref{eq:ASW400}) and a fifth missing identity which was originally conjectured by Feigin--Foda--Welsh \cite{FFW08}. All these modulo 7 identities have manifestly positive sum sides. \cite{CW19} has been the catalyst for the recent developments. Ablinger and the author \cite{qFunctions} implemented the Corteel--Welsh's cylindric partitions related functional equations in their symbolic computation implementation \texttt{qFunctions} to be able to exploit this combinatorial idea using formal manipulation and computer algebra techniques. Corteel, Dousse and the author \cite{CDU} later proved the modulo 8 identities that arise from the cylindric partitions paradigm with the help of this implementation. One of such identities is as follows (see Theorem~1.6 in \cite{CDU}).

\begin{theorem}[Corteel--Dousse--U., 2021]
\label{th:CDU500}
\begin{equation}\label{eq:CDU500}\sum_{\substack{r_1\geq s_1\geq r_2\geq 0\\r_1\geq s_2\geq 0}} \frac{q^{r_1^2 -r_1 s_1 + s_1^2 + r_2^2 + s_2^2  +s_1 s_2 + r_1 + r_2 + s_1 + s_2}}{(q;q)_{r_1}}{r_1\brack s_1}_q{r_1\brack s_2}_q{s_1\brack r_2}_q =\frac{1}{\theta(q^2,q^3,q^3,q^4,q^4,q^5;q^{10})}\end{equation}
\end{theorem}

Unlike \eqref{eq:ASW500}, \eqref{eq:CDU500} has a manifestly positive sum-side. Shortly after \cite{CDU}, in late 2021, Warnaar \cite{War21} come up with many beautiful conjectures for manifestly positive sum-sides related to higher moduli (not divisible by 3). In 2022, Tsuchioka \cite{Tsu22} proved manifestly positive sum-sides for modulus $6$ using finite-automata and automated proofs. He was also able to analyze the structure of relevant level 3 standard modules for the affine Lie algebra $A_2^{(1)}$.

In a different vein, Bridges and the author studied weighted versions of cylindric partitions as well as cylindric partitions into distinct parts in \cite{BU}.

Earlier in 2022, Kanade and Russell \cite{KR} aimed (and succeeded) at conjecturing $A_2$ Rogers--Ramanujan type identities in the form of Andrews--Schilling--Warnaar instead of aiming for manifestly positive sum-sides. They were able to make explicit claims for each modulus $\geq 5$. They proved the cases for moduli $5,6,7,8$ and 10. Their exploration came to an end due to the increasing computational difficulties.

In this paper, we approach the conjectures of Kanade--Russell by changing the computational techniques used. We prove all modulo 11 and 13 $A_2$ Rogers--Ramanujan identities coming from the cylindric partitions paradigm. Two such identities are as follows:

\begin{theorem}\label{th:mod11ex}
\begin{align*}\sum_{\substack{r_1\geq r_2\geq r_3\geq 0\\s_1\geq s_2\geq s_3\geq0 }} &\frac{q^{r_1^2 - r_1 s_1+ s_1^2+r_2^2 - r_2 s_2+ s_2^2 +r_3^2 + r_3 s_3+ s_3^2+r_1+r_2+r_3+s_1+s_2+s_3}}{ (q;q)_{r_1-r_{2}}(q;q)_{r_2-r_{3}}(q;q)_{r_{3}}(q;q)_{s_1-s_{2}}(q;q)_{s_2-s_{3}}(q;q)_{s_{3}}(q;q)_{r_{3}+s_{3}+1}}\\&\hspace{7cm} = \frac{1}{(q;q)_\infty} \frac{1}{\theta(q^{2},q^{3},q^{3},q^{4},q^{4},q^{5},q^{5};q^{11})}.\end{align*}
\end{theorem}

\begin{theorem}\label{th:mod13ex}
\begin{align*}\sum_{\substack{r_1\geq r_2\geq r_3\geq 0\\s_1\geq s_2\geq s_3\geq0 }} &\frac{q^{r_1^2 - r_1 s_1+ s_1^2+r_2^2 - r_2 s_2+ s_2^2 +r_3^2 - r_3 s_3+ s_3^2+r_1+r_2+r_3+s_1+s_2+s_3}}{ (q;q)_{r_1-r_{2}}(q;q)_{r_2-r_{3}}(q;q)_{r_{3}}(q;q)_{s_1-s_{2}}(q;q)_{s_2-s_{3}}(q;q)_{s_{3}}(q;q)_{r_{3}+s_{3}+1}}\\&\hspace{7cm} = \frac{1}{(q;q)_\infty} \frac{1}{\theta(q^{2},q^{3},q^{3},q^{4},q^{4},q^{5},q^{5},q^6,q^6;q^{13})}.\end{align*}
\end{theorem}

The organization of this paper is as follows. In Section~\ref{sec:notation}, we introduce cylindric partitions, the relevant results and the conjectures of Kanade--Russell of which we prove some cases of. Section~\ref{sec:proofmethod} is dedicated to rewording the conjectures and the description ot the proof methodology. In Sections~\ref{sec:M11} and \ref{sec:M13} we present the proofs of the modulo 11 and 13 $A_2$ Rogers--Ramanujan identities in Andrews--Schilling--Warnaar form, respectively. We outline some natural questions and mathematical challenges that arise from this work in Section~\ref{sec:future}. Section~\ref{sec:computations} is reserved for a discussion on how the computerized proofs have been carried in earlier work \cite{CDU, KR} and this paper and what future improvements can be done to take us further mathematically.

\section*{Acknowledgement}

The author would like to thank the workshop on cylindric partitions group that came together in November 2022 in Linz for all the stimulating discussions. In particular, the author would like to thank Shashank Kanade for suggesting that the researchers working on cylindric partitions should come together and join forces in the first place, and for all his comments on this manuscript. The author would also like to thank Christian Koutschan his encouragement of the author in the necessary implementations.

Research of the author is partly supported by EPSRC grant number EP/T015713/1 and partly by FWF grant P-34501N.

\section{Necessary definitions}\label{sec:notation}

We shall start with the definition of a cylindric partition.

\begin{definition}
 A \textit{cylindric partition} is made up of a composition $c=(c_1,c_2,\dots,c_r)$ called \textit{profile} with $r$ parts, and a vector $\pi = (\pi^{(1)}, \pi^{(1)}, \dots, \pi^{(r)})$ consisting of $r$ partitions $\pi^{(i)} = (\pi^{(i)}_1, \pi^{(i)}_2,\dots)$, that satisfy the inequalities \[\pi^{(i)}_j\geq \pi^{(i+1)}_{j+c_{i+1}}\quad \text{and}\quad \pi^{(r)}_j\geq \pi^{(1)}_{j+c_1}.\]
 \end{definition}
 
For example, the vector partition $\pi=\{(1,1,1,1), (4,3,1), (2,2)\}$ together with the profile $(2,0,2)$ is a cylindric partition. Note that the same vector partition can also satisfy the cylindric partition inequalities with different profiles. For example, $\pi$ is also a cylindric partition for profiles $(2,0,0)$, $(2,0,1)$, etc. We can define the total size of a cylindric partition $\pi$ as the sum of all the sizes of the partitions included. We denote the total size, once again, by $|\pi|$. Let $c$ be a composition and let $\P_c$ be the set of all vector partitions that are cylindric partitions with profile $c$. 

For a given profile $c$, let $\mathcal{P}_c$ be the set of all cylindric partitions with profile $c$. Let \[F_{c}(z,q) := \sum_{\pi\in \P_c} z^{\max(\pi)}q^{|\pi|},\] the bivariate generating function for the number of cylindric partitions where the exponents of $z$ and $q$ are keeping record of the largest parts size and the total of the parts in $\pi$, respectively. Borodin \cite{Bor07} showed that when $z=1$, $F_c(z,q)$ generating functions have product formula.

\begin{theorem}[Borodin, 2007]
\label{th:Borodin}
Let $r$ and $l$ be positive integers, and let $c=(c_1,c_2,\dots,c_r)$ be a composition of $l$. Define $m:=r+l$ and $s(i,j) := c_i+c_{i+1}+\dots+ c_j$. Then,
\begin{equation}
\label{BorodinProd}
F_c(1,q) = \frac{1}{(q^m;q^m)_\infty} \prod_{i=1}^r \prod_{j=i}^r \prod_{k=1}^{c_i} \frac{1}{(q^{k+j-i+s(i+1,j)};q^m)_\infty} \prod_{i=2}^r \prod_{j=2}^i \prod_{k=1}^{c_i} \frac{1}{(q^{m-k+j-i-s(j,i-1)};q^m)_\infty}.
\end{equation}
\end{theorem}

Focusing on replacing the largest part in a given cylindric partition, Corteel--Welsh \cite{CW19} defined a $q$-difference equation for $F_c(z,q)$. This functional equation relates $F_c(z,q)$ with other generating functions $F_{c^*}(z,q)$ where $\#(c)=\#(c^*)$ and $|c|=|c^*|$. Let $c=(c_1, \dots , c_r)$ (with the convention that $c_0=c_r$) be a given composition and define $I_c$ to be the set of indices for the non-zero entries in $c$. Given a non-empty subset $J\subseteq I_c$, the composition $c(J) = (c_1(J), \dots , c_r(J))$ is defined by:
\begin{equation}
\label{eq:c(J)}
c_i(J):= \begin{cases}
c_i-1 &\text{if $i \in J$ and $i-1 \notin J$},\\
c_i+1 &\text{if $i \notin J$ and $i-1 \in J$},\\
c_i &\text{otherwise}.
\end{cases}
\end{equation}
Then the explicit $q$-difference equation $F_c(z,q)$ satisfies is as follows.

\begin{theorem}[Corteel--Welsh, 2019] 
\label{th:CW}
For any profile $c$, 
\begin{equation}
\label{CorteelRec}
F_c(z,q) = \sum_{\emptyset\subset J\subseteq I_c} (-1)^{|J|-1} \frac{F_{c(J)}(z q^{|J|},q)}{(1-z q^{|J|})},
\end{equation}
 with the initial conditions $F_c(0,q)=F_c(z,0)=1$.
\end{theorem}

Let $c$ be a profile and $c'$ be a cyclic shift of $c$. There is a clear one-to-one correspondence between cylindric partitions in $\mathcal{P}_c$ and $\mathcal{P}_{c'}$ by cyclically shifting the vector of partitions counted in $\mathcal{P}_{c}$. This is enough to see that the generating functions for these sets of cylindric partitions are equal, i.e. $F_c(z,q) = F_{c'}(z,q)$. Therefore, we can cyclically shift the profiles and lower the number of  (seemingly different) generating functions that appear in the coupled system of $q$-difference equations.

We can also normalize \eqref{CorteelRec} and get an equivalent $q$-difference equation. For example, let
$$G_c(z,q):= (zq;q)_{\infty} F_c(z,q).$$ 
The equation \eqref{CorteelRec} is equivalent to

 \begin{equation}\label{eq:Grec}
G_c(z,q) = \sum_{\emptyset\subset J\subseteq I_c} (-1)^{|J|-1} (zq;q)_{|J|-1} G_{c(J)}(z q^{|J|},q),\end{equation} with the initial conditions $G_c(0,q)=G_c(z,0)=1$. This $q$-difference equation \eqref{eq:Grec} with polynomial coefficients, in practice, played a central role in the proofs of modulo 7 and modulo 8 cylindric partition with 3-part profile identities in \cite{CW19} and \cite{CDU}, respectively. Weighted versions of \eqref{CorteelRec} and \eqref{eq:Grec} are later presented in \cite{BU}. In \cite{KR}, Kanade--Russell decided to change the initial conditions of \eqref{eq:Grec} slightly. While this does not change the $q$-difference equations, this lead to the conjectural discovery of explicit formulas for most of these 3-part profile cylindric partition generating functions. Let \begin{equation}\label{eq:HtoF} H_c(z,q):= \frac{(zq;q)_{\infty}}{(q;q)_\infty} F_c(z,q).\end{equation} Then $H_c(z,q)$ satisfies the same $q$-difference equation as $G_c(z,q)$, namely \begin{equation}\label{eq:Hrec}
H_c(z,q) = \sum_{\emptyset\subset J\subseteq I_c} (-1)^{|J|-1} (zq;q)_{|J|-1} H_{c(J)}(z q^{|J|},q),\end{equation} with the initial conditions $H_c(0,q)=1/(q;q)_\infty$ and $H_c(z,0)=1$. 

From this point forward we only focus on cylindric partition profiles with 3-parts.

Let $k\geq 2$, let \begin{equation}\label{eq:rho_sigma}\rho = (\rho_1,\rho_2,\dots, \rho_{k-1})\in\mathbb{Z}^{k-1}, \ \text{and}\ \sigma = (\sigma_1,\sigma_2,\dots, \sigma_{k-1})\in\mathbb{Z}^{k-1}\end{equation} and define 
\begin{align}
  \label{eq:Sm1}  S_{3k-1}(z;\rho | \sigma) &= \sum_{r,s \in \mathbb{Z}^{k-1}_{\geq 0}} z^{r_1} \frac{q^{ \sum_{i=1}^{k-1} r_i^2 - r_i s_i +s_i^2 + \rho_i r_i + \sigma_i s_i}}{\prod_{i=1}^{k-2} (q;q)_{r_i-r_{i+1}}(q;q)_{s_i-s_{i+1}}} \frac{q^{2r_{k-1}s_{k-1}}}{(q;q)_{r_{k-1}}(q;q)_{s_{k-1}}(q;q)_{r_{k-1}+s_{k-1}+1}},\\
\label{eq:S0}    S_{3k}(z;\rho | \sigma) &= \sum_{r,s \in \mathbb{Z}^{k-1}_{\geq 0}} z^{r_1} \frac{q^{\sum_{i=1}^{k-1} r_i^2 - r_i s_i+ s_i^2 + \rho_i r_i + \sigma_i s_i}}{\prod_{i=1}^{k-2} (q;q)_{r_i-r_{i+1}}(q;q)_{s_i-s_{i+1}}} \frac{1}{(q;q)_{r_{k-1}+s_{k-1}}(q;q)_{r_{k-1}+s_{k-1}+1}} {r_{k-1}+s_{k-1}\brack r_{k-1}}_{q^3},\\
\label{eq:Sp1}    S_{3k+1}(z;\rho | \sigma) &= \sum_{r,s \in \mathbb{Z}^{k-1}_{\geq 0}} z^{r_1} \frac{q^{\sum_{i=1}^{k-1} r_i^2 - r_i s_i+ s_i^2 + \rho_i r_i + \sigma_i s_i}}{\prod_{i=1}^{k-2} (q;q)_{r_i-r_{i+1}}(q;q)_{s_i-s_{i+1}}} \frac{1}{(q;q)_{r_{k-1}}(q;q)_{s_{k-1}}(q;q)_{r_{k-1}+s_{k-1}+1}}.
\end{align}

Let \begin{equation}\label{eq:ei} e_i = (\underbrace{0,0,\dots,0}_{i},1,1,\dots,1)\in \Z^{k-1}\quad\text{and}\quad \delta_i := (\delta_{ij})_{1\leq j\leq k-1} \in \mathbb{Z}^{k-1}.\end{equation} It is easy to see that \begin{equation}\label{eq:HzTransform} S_m(zq^n; \rho | \sigma) = S_m(z; \rho + n \delta_1 | \sigma)\end{equation} for any $n\in \mathbb{Z}$.

Kanade--Russell conjectured that for any fixed $k\geq 3$, the $H_{(c_1,c_2,c_3)}(z,q)$ can be expressed as linear combinations of the $S_{m}(z;\rho |\sigma)$ functions. Precisely they claimed the following.

\begin{conjecture}[Kanade--Russell, 2022]\label{conj:Hconj} Let $k\geq 3$ and $|c| + 3 = m = 3k+\{-1,0,1\}$. Using cyclis symmetries asusme that $c_1\geq c_2, c_3$. If $c_2, c_3 \leq k-1$, then
\begin{equation}\label{eq:Hconj}
    H_{(c_1,c_2,c_3)} (z,q) = \left\{\begin{array}{ll}
         S_m( z; e_{c_2} |e_{c_3})-qS_m( z; e_{c_2-1} |e_{c_3-1}), & c_2,c_3>0,  \\
         S_m( z;e_{c_2} | e_{0}), & c_3=0, \\
         S_m( z; e_{0} |e_{c_3})-q(1-z)S_m( z; e_{0}+\delta_0 |e_{c_3-1}), & c_2=0, c_3\not=0,
    \end{array}  \right. 
\end{equation} where $e_i$ and $\delta_{i}$ be defined as in \eqref{eq:ei}.
\end{conjecture}

The explicit claims that Conjecture~\ref{conj:Hconj} provide do not cover all the functions $H_c(z,q)$ with $|c|+3=m$. It does provide enough claims to recover explicit expression claims. How to find the conjectural $S_m( z;\rho|\sigma)$ equivalents of the other $H_c(z,q)$ that appear in the coupled $q$-difference equation system is explained in \cite{KR}. The profiles related to the functions to be recovered are called "\textit{under-the-line}" profiles by Kanade--Russell. We will also call these profiles as such while we ignore to explain anything about the line. These under-the-line profile related functions can have shifts of $z$ in the $S_m( z;\rho|\sigma)$ language. These shifts are inherited from the $q$-difference equations \eqref{eq:Hrec}. One can use \eqref{eq:HzTransform} to clear all the shifts in $z$. Therefore, from now on in all our expressions we will translate any $z$ shift of $S_m(z;\rho|\sigma)$ using \eqref{eq:HzTransform} and this way ignore any and all shifts in $z$. To further emphasize this moving forward on we suppress the variable $z$ from our notation and write \[ S_m(\rho|\sigma):=  S_m(z;\rho|\sigma).\]
 
Proof of Conjecture~\ref{conj:Hconj} (and its extension to all 3-part profiles with total $m-3$) requires one to show that the expressions in $S_m(\rho|\sigma)$ are the correct expressions for the respective $H_c(z;q)$ functions. This can be done by showing that the expressions in $S_m(\rho|\sigma)$ satisfies the same recurrence relation specified by \eqref{eq:Hrec} and the initial conditions of the expression holds. In \cite{KR}, it is already proven that for $c_2, c_3 \leq k-1$, the conjectural formulas of \eqref{eq:Hconj} all satisfy the necessary initial conditions \begin{equation}
    \label{eq:inits}H_{c}(z,0)=1\quad\text{and}\quad H_c(0,z) = 1/(q;q)_\infty.
\end{equation}

It was noted in the \cite[Lemma 9.1, Lemma 9.2]{KR} that $S_m(\rho|\sigma)$ functions satisfy the following list of recurrences.

\begin{lemma}[Kanade--Russell, 2022]\label{lemma:recs} Let $k\geq 3$, let $m = 3k+\{-1,0,1\}$ and let $\delta_i := (\delta_{ij})_{1\leq j\leq k-1} \in \mathbb{Z}^{k-1}$, where $\delta_{ij}$ is the Kronecker delta function. The following recurrence relations follow for all $1\leq i\leq k-2$,
\begin{align}
\label{R1}\tag{$R_1^{(i)}(\rho|\sigma)$}  S_m(\rho | \sigma) &- S_m(\rho +\delta_i-\delta_{i+1} | \sigma) - zq^{i+\sum_{j=1}^i \rho_j}S_m(\rho + 2 \sum_{j=1}^i \delta_j | \sigma -\sum_{j=1}^i \delta_j)=0,\\
\label{R2}\tag{$R_2^{(i)}(\rho|\sigma)$} S_m(\rho | \sigma) &- S_m(\rho  | \sigma+\delta_i-\delta_{i+1}) - zq^{i+\sum_{j=1}^i \sigma_j}S_m(\rho -\sum_{j=1}^i \delta_j | \sigma+ 2 \sum_{j=1}^i \delta_j )=0.
\end{align}
\begin{enumerate}[i.]
    \item If $m \equiv -1 \pmod 3$, \begin{enumerate} [a)]
        \item and if $\sigma_{k-1} = 0$, then \begin{equation}\label{R3m1}\tag{$R_3(\rho|\sigma)$}S_m(\rho | \sigma) - S_m(\rho | \sigma + \delta_{k-1}) - q S_m(\rho +\delta_{k-1} | \sigma + \delta_{k-1}) +q S_m(\rho +\delta_{k-1} | \sigma+\delta_{k-2} + \delta_{k-1})=0.\end{equation}
        \item and if $\rho_{k-1} = 0$, then \begin{equation}\label{R4m1}\tag{$R_4(\rho|\sigma)$}S_m(\rho | \sigma) - S_m(\rho + \delta_{k-1} | \sigma ) - q S_m(\rho +\delta_{k-1} | \sigma + \delta_{k-1}) +q S_m(\rho+\delta_{k-2} +\delta_{k-1} | \sigma + \delta_{k-1})=0.\end{equation}
    \end{enumerate}
    \item If $m \equiv 0 \pmod 3$, then \begin{align}
    \label{R30}\tag{$R_3(\rho|\sigma)$}S_m(\rho | \sigma) &- (1+q) S_m(\rho + \delta_{k-1}| \sigma + \delta_{k-1}) + q S_m(\rho +2\delta_{k-1} | \sigma + 2\delta_{k-1})\\ \nonumber &-z q^{k-1+\sum_{j=1}^{k-1} \rho_j} S_m(\rho +2\sum_{j=1}^{k-1} \delta_j | \sigma-\sum_{j=1}^{k-1} \delta_j )\\ \nonumber &- q^{k-1+\sum_{j=1}^{k-1} \sigma_j} S_m(\rho -\sum_{j=1}^{k-2} \delta_j  +2 \delta_{k-1}| \sigma+2\sum_{j=1}^{k-1} \delta_j )=0.\\
    \label{R40}\tag{$R_4(\rho|\sigma)$}S_m(\rho | \sigma) &- (1+q) S_m(\rho + \delta_{k-1}| \sigma + \delta_{k-1}) + q S_m(\rho +2\delta_{k-1} | \sigma + 2\delta_{k-1})\\ \nonumber &-z q^{k-1+\sum_{j=1}^{k-1} \rho_j} S_m(\rho +2\sum_{j=1}^{k-1} \delta_j | \sigma-\sum_{j=1}^{k-2} \delta_j +2\delta_{k-1})\\ \nonumber &- q^{k-1+\sum_{j=1}^{k-1} \sigma_j} S_m(\rho -\sum_{j=1}^{k-1} \delta_j  | \sigma+2\sum_{j=1}^{k-1} \delta_j )=0.
    \end{align}
    \item If $m \equiv 1 \pmod 3$, then
    \begin{align}
    \label{R31}\tag{$R_3(\rho|\sigma)$}  S_m(\rho | \sigma) &- S_m(\rho  | \sigma+ \delta_{k-1}) - q S_m(\rho + \delta_{k-1} | \sigma + 2\delta_{k-1})\\\nonumber&+q S_m(\rho + \delta_{k-1} | \sigma + 2\delta_{k-1}) -q^{k-1+\sum_{j=1}^{k-1} \sigma_j} S_m(\rho -\sum_{j=1}^{k-1} \delta_j  | \sigma+2\sum_{j=1}^{k-1} \delta_j )=0 \\
        \label{R41}\tag{$R_4(\rho|\sigma)$}  S_m(\rho  | \sigma) &- S_m(\rho + \delta_{k-1} | \sigma) - q S_m(\rho + \delta_{k-1} | \sigma + \delta_{k-1})\\\nonumber&+q S_m(\rho +2 \delta_{k-1} | \sigma + \delta_{k-1}) -zq^{k-1+\sum_{j=1}^{k-1} \rho_j} S_m(\rho +2\sum_{j=1}^{k-1} \delta_j  | \sigma-\sum_{j=1}^{k-1} \delta_j )=0 .
    \end{align}
\end{enumerate}
\end{lemma}

Then they made the following claim (see \cite[Conjecture~9.3]{KR}).

\begin{conjecture}[Kanade--Russell, 2022]\label{KRconjecture}
In each modulus $m\geq 5$, the relations \eqref{R1}-\eqref{R40} are enough to prove recurrences necessary for the proof of Conjecture~\ref{conj:Hconj}.
\end{conjecture}

We find this conjecture highly sensible. For all $m\geq 5$, the explicit $S_m$'s are $2\lfloor m/3 \rfloor$-fold sums. Same is true for the number of distinct functional equations \eqref{R1}-\eqref{R40}. One can easily check that these relations are distinct by comparing the first two terms in each left-hand side. Each second term corresponds to a canonical shift in one of the summation variables. One would expect to see every relation that the $S_m(\rho|\sigma)$ functions satisfy to be translated and recovered as a combination the relations \eqref{R1}-\eqref{R40}. Hence, if the claims of Conjecture~\ref{conj:Hconj} are correct, for any fixed profile $c$ the coupled $q$-difference equations \eqref{eq:Hrec} written using the explicit claims of \eqref{eq:Hconj} (together with the ``\textit{under-the-line}" expressions) can be recovered as a combination of the relations \eqref{R1}-\eqref{R40}. 

\section{Proof Methodology}\label{sec:proofmethod}

Conjecture~\ref{KRconjecture} can be rephrased as a set inclusion question. Let $m = 3k +\{-1,0,1\}$ with $k\geq 3$, $\rho$ and $\sigma$ as in \eqref{eq:rho_sigma} and $1\leq i\leq k-2$. Define \begin{equation}\label{eq:Im}I_{m} := \langle R^{(i)}_1(\rho|\sigma),R^{(i)}_2(\rho|\sigma), R_3(\rho|\sigma),R_4(\rho|\sigma)\rangle,\end{equation}  the ideal generated by the left-hand sides of the recurrences \eqref{R1}-\eqref{R41} as polynomials in the ring $\mathbb{Z}(\!(q, z)\!)[S_m(\rho|\sigma)]$. Here which $R_3(\rho|\sigma)$ and $R_4(\rho|\sigma)$ to be included in $I_{m}$ is to be understood by the residue class of $m$ modulo 3. Recall that $\rho$ and $\sigma$ are integer vectors with $k-1$ entries. Therefore the ring $\mathbb{Z}(\!(q, z)\!)[S_m(\rho|\sigma)]$ is a formal polynomial ring defined on a countable set of variables.

For any given fixed $m = 3k +\{-1,0,1\}$ with $k\geq 3$, let the set of all the coupled system of $q$-difference equations \eqref{eq:Hrec} for the profiles $c$ with $|c|+3=m$ be $\mathcal{H}_m$. Any relation in $\mathcal{H}_m$ can be written in $S_m(\rho|\sigma)$ functions using the Conjecture~\ref{conj:Hconj} (and the paragraph below it). Let $\mathcal{S}_m$ be the set of all relations in $\mathcal{H}_m$ written in the conjectural $S_m(\rho|\sigma)$ form. 

Now we can write Conjecture~\ref{KRconjecture} in its equivalent form:

\begin{conjecture}\label{conj:Ideal} Let $m = 3k +\{-1,0,1\}$ with $k\geq 3$ be fixed. \[\forall h\in \mathcal{S}_m,\text{ we have}\ \ h \in I_{m}.\]
\end{conjecture}

The infinite set $\{ R^{(i)}_1(\rho|\sigma),R^{(i)}_2(\rho|\sigma), R_3(\rho|\sigma),R_4(\rho|\sigma)\, :\, 1\leq i\leq k-2, \rho,\sigma\in\mathbb{Z}^{k-1}\}$ that spans $I_m$ has non-trivial relations within itself and not all the elements of this set are generators of $I_{m}$. However, we do not know an exact pattern of which elements are related at the moment. Nevertheless, it is easy to understand that $I_m$ is generated by infinitely many elements since $\rho$ and $\sigma\in\Z^{k-2}$. 

On the other hand, for any fixed $m$, the $S_m(\rho|\sigma)$ functions that appear within the formulas from $\mathcal{S}_m$ make up a finite list. One can easily find explicit bounds for the entries of vectors $\rho$ and $\sigma$ such that every $S_m(\rho|\sigma)$ that appear in $\mathcal{S}_m$ is within the bounds. This observation suggests that instead of attempting to prove Conjecture~\ref{conj:Ideal}, we can instead go after a stronger conjecture that is more suitable for computations. To that end, let $[N]:=\{-N,\dots,-1,0,1,\dots,N\}$ and we define \[I_{m,N} := \langle \{ R^{(i)}_1(\rho|\sigma),R^{(i)}_2(\rho|\sigma), R_3(\rho|\sigma),R_4(\rho|\sigma)\, :\, \rho,\sigma \in [N]^{k-1}   \} \rangle \subset I_m.\]

With this definition we form the stronger conjecture:

\begin{conjecture}\label{conj:FinIdeal} Let $m = 3k +\{-1,0,1\}$ with $k\geq 3$ be fixed. There is some $N\in \mathbb{N}$ such that \[\forall h\in \mathcal{S}_m,\text{ we have}\ \ h \in I_{m,N}.\]
\end{conjecture}

Since $I_{m,N}\subset I_m$, it is clear that Conjecture~\ref{conj:FinIdeal} implies Conjecture~\ref{conj:Ideal}.

Finally we transferred the open problems into a linear algebra setting, and we can approach it as such. 

Let $m$ and $N$ be fixed. we can order all the $S_m(\rho|\sigma)$ that appears in the spanning set of $I_{m,N}$ and write in a column vector $\vec{s}$. Then the matrix $\mathbf{A}$ is uniquely defined by \[ \mathbf{A} \vec{s} = \vec{0}_\mathbf{A},\] where $\vec{0}_\mathbf{A}$ is the colum vector with the same number of rows as $\mathbf{A}$. Every row of $\mathbf{A}$, corresponds to a functional relation $R_j(\sigma|\rho) \in  \{ R^{(i)}_1(\rho|\sigma),R^{(i)}_2(\rho|\sigma), R_3(\rho|\sigma),R_4(\rho|\sigma)\, :\, \rho,\sigma \in [N]^{k-1}  \}$ and every column of $\mathbf{A}$ corresponds to the coefficients of $S_m(\rho|\sigma)$. Also observe that $\mathbf{A}$ is a finite dimensional matrix with entries in $\mathbb{Z}[q,z]$.

One can use Gaussian elimination on $\mathbf{A}$. Any non-trivial relation within the functional relations \eqref{R1}-\eqref{R41} within the defining bounds of $\mathbf{A}$ would yield 0 rows. Let $\mathbf{B}$ be the matrix consisting of non-zero rows of $\mathbf{A}$ after the Gaussian elimination is performed. It should still be clear that \[\mathbf{B}\vec{s} = \vec{0}_\mathbf{B}.\] Moreover, the ideal $I_{m,N}$ is generated by the equations that appear in $ \mathbf{B}\vec{s} = \vec{0}$.

Therefore, for any element of $h\in\mathcal{S}_m$ one can check whether that element is in $I_{m,N}$ by simply writing that relation as a row vector $\vec{h}$ (with respect to the vector $\vec{s}$, i.e. $\\vec{h}$ is defined by $h := [\vec{h}\vec{s}=0]$), add the row vector $\vec{h}$ to $\mathbf{B}$ and perform Gaussian elimination to this new matrix. If the Gaussian elimination yields a zero row, this means that $\vec{h}$ is a linear combination of rows in $\mathbf{B}$, or equivalently this means $h\in I_{m,N}$. If no zero row appears, then $h\not\in I_{m,N}$.

This approach is clearly algorithmic. Furthermore, termination of the algorithm and a definitive answer among the termination are both guaranteed. Top it all up, the explicit combination of \eqref{R1}-\eqref{R41} functional equations that is equivalent to a given  $h\in\mathcal{S}_m$ is also easy to find. One only needs to use an augmented version of $\mathbf{A}$ where one more column is added to keep track of the name of the relations \eqref{R1}-\eqref{R41} while doing the row reductions.

After these considerations, proof of Conjectures~\ref{conj:FinIdeal} (and consequently Conjectures~\ref{conj:Hconj}, \ref{KRconjecture}, and \ref{conj:Ideal}) comes down to experimentally identifying an $N$ and being able to perform the Gaussian elimination calculations.

\section{Modulo 11 Identities}\label{sec:M11}

Let $m=11$ ($= 3 k -1$ with $k=4$), for this family of identities $\rho$ and $\sigma \in \mathbb{Z}^3$. There are a total of 15 essentially unique 3 part compositions of 8 that appear in the coupled $q$-difference system \eqref{eq:Hrec}. Conjecture~\ref{conj:Hconj} suggests that the following sum representations for $H_c(z,q)$ hold for all but one of these:

\begin{equation}\label{eq:mod11list}
\begin{array}{ll}
H_{(8, 0, 0)}(z,q) &= S_{11} ( (1, 1, 1)\, |\, (1, 1, 1)),\\
H_{(7, 1, 0)}(z,q) &= S_{11} ( (0, 1, 1)\, |\, (1, 1, 1)),\\ 
H_{(7, 0, 1)}(z,q) &= 
 S_{11} ( (1, 1, 1)\, |\, (0, 1, 1)) - q (1 - z) S_{11} ( (2, 1, 1)\, |\, (1, 1, 1)),\\
H_{(6, 2, 0)}(z,q) &= S_{11} ( (0, 0, 1)\, |\, (1, 1, 1)),\\ 
H_{(6, 1, 1)}(z,q) &= 
 S_{11} ( (0, 1, 1)\, |\, (0, 1, 1)) - q S_{11} ( (1, 1, 1)\, |\, (1, 1, 1)),\\ 
H_{(6, 0, 2)}(z,q) &= 
 S_{11} ( (1, 1, 1)\, |\, (0, 0, 1)) - q (1 - z) S_{11} ( (2, 1, 1)\, |\, (0, 1, 1)),\\  
H_{(5, 3, 0)}(z,q) &= S_{11} ( (0, 0, 0)\, |\, (1, 1, 1)),\\ 
H_{(5, 2, 1)}(z,q) &= 
 S_{11} ( (0, 0, 1)\, |\, (0, 1, 1)) - q S_{11} ( (0, 1, 1)\, |\, (1, 1, 1)),\\
H_{(5, 1, 2)}(z,q) &= 
 S_{11} ( (0, 1, 1)\, |\, (0, 0, 1)) - q S_{11} ( (1, 1, 1)\, |\, (0, 1, 1)),\\
H_{(5, 0, 3)}(z,q) &= 
 S_{11} ( (1, 1, 1)\, |\, (0, 0, 0)) - q (1 - z) S_{11} ( (2, 1, 1)\, |\, (0, 0, 1)),\\ 
H_{(4, 3, 1)}(z,q) &= 
 S_{11} ( (0, 0, 0)\, |\, (0, 1, 1)) - q S_{11} ( (0, 0, 1)\, |\, (1, 1, 1)),\\ 
H_{(4, 2, 2)}(z,q) &= 
 S_{11} ( (0, 0, 1)\, |\, (0, 0, 1)) - q S_{11} ( (0, 1, 1)\, |\, (0, 1, 1)),\\ 
H_{(4, 1, 3)}(z,q) &= 
 S_{11} ( (0, 1, 1)\, |\, (0, 0, 0)) - q S_{11} ( (1, 1, 1)\, |\, (0, 0, 1)),\\ 
H_{(3, 3, 2)}(z,q) &= 
 S_{11} ( (0, 0, 0)\, |\, (0, 0, 1)) - q S_{11} ( (0, 0, 1)\, |\, (0, 1, 1)).\\ 
\end{array}
\end{equation}

Only $H_{(4,4,0)}(z,q)$ misses a claimed formula and that can be recovered by the $q$-difference equations \eqref{eq:Hrec}. We know that $H_{(4,4,0)}(z,q)$ satisfies
\begin{equation}\label{eq:H440}H_{(4, 4, 0)}(z,q) +(1 - q z) H_{(4, 1, 3)}(q^2 z,q) - H_{(4, 3, 1)}(q z,q)  -  H_{(5, 0, 3)}(q z,q) = 0.\end{equation}
Using the conjectured series equivalents \eqref{eq:mod11list} of $H_{(4, 1, 3)}(z,q)$, $H_{(4, 3, 1)}(z,q)$ and $H_{(5,0, 3)}(z,q)$, we see that
\begin{align}\nonumber H_{(4,4,0)}(z,q) &= -(S_{11}(qz;(0, 0, 0)|(0, 1, 1)) - q S_{11}(qz;(0, 0, 1)|(1, 1, 1))) \\\label{eq:H440z}&+ (1 -q z) (S_{11}(q^2 z;(0, 1, 1)|(0, 0, 0)) - q S_{11}(q^2 z;(1, 1, 1)|(0, 0, 1))) \\ \nonumber &- (S_{11}(q z; (1, 1, 1)|(0, 0, 0)) - q (1 - z) S_{11}(q z(2, 1, 1)|(0, 0, 1))).
\end{align} Notice that we used the shifts in the variable $z$ in \eqref{eq:H440z}. We clear these shifts by employing \eqref{eq:HzTransform}. This yields an explicit claim for $H_{(4,4,0)}(z,q)$: \begin{equation}\label{eq:H440exp}H_{(4,4,0)}(z,q) =S_{11}((1, 0, 0)\, |\, (0, 1, 1)) - q S_{11}((1, 0, 1)\, |\, (1, 1, 1)) + q z S_{11}(2, 1, 1)\, |\, (0, 0, 0)), \end{equation} with no shifts in $z$, where the $S_{11}(\rho|\sigma)$ functions fit the forms in Lemma~\ref{lemma:recs}. 

We can also see that $H_{(4,4,0)}(z,q)$ satisfies the necessary initial conditions \eqref{eq:inits}. The initial condition $H_{(4,4,0)}(z,0)=1$ is immediate by \eqref{eq:H440exp} and \eqref{eq:Sm1}. We can also see that $H_{(4,4,0)}(0,q) = 1/(q;q)_\infty$ by plugging in $z=0$ in \eqref{eq:H440} and using the initial conditions of the other proven initial conditions \eqref{eq:inits} for the functions in \eqref{eq:H440}.

Our proof routine explained in Section~\ref{sec:proofmethod} can start once all the normalized generating functions $H_c(z,q)$'s are (conjecturally) translated in the $S_{11}((a_1,a_2,a_3)|(b_1,b_2,b_3))$ language. It is easy to see that the following four recurrences,
\begin{align*}
H_{(8, 0, 0)}(z,q)& -H_{(7, 1, 0)}(q z,q) = 0,\\
H_{(7, 0, 1)}(z,q)&  -H_{(6, 1, 1)}(q z,q) + (1 - q z) H_{(7, 1, 0)}(q^2 z,q) -
  H_{(8, 0, 0)}(q z,q) = 0,\\
H_{(6, 0, 2)}(z,q)& -H_{(5, 1, 2)}(q z,q) + (1 - q z) H_{(6, 1, 1)}(q^2 z,q) -
  H_{(7, 0, 1)}(q z,q) =0,\\
H_{(5, 0, 3)}(z,q)&-H_{(4, 1, 3)}(q z,q)  + (1 - q z) H_{(5, 1, 2)}(q^2 z,q) -
  H_{(6, 0, 2)}(q z,q) = 0,\\
\end{align*} 
trivializes to $0=0$ once the terms on the left-hand sides are written in $S_{11}((a_1,a_2,a_3)|(b_1,b_2,b_3))$ using \eqref{eq:mod11list} and \eqref{eq:HzTransform}. Therefore, these relations are trivially in $I_{11}$, the ideal generated by the functional relations of the $S_{11}((a_1,a_2,a_3)|(b_1,b_2,b_3))$ series.

Recall that we used the coupled $q$-difference equation \eqref{eq:H440} to make an explicit claim for $H_{(4,4,0)}(z,q)$. Hence, the functional relation of $H_{(4,4,0)}(z,q)$ also trivializes to $0=0$ once written in the claimed $S_{11}(\rho|\sigma)$ forms.  The very claim \eqref{eq:H440exp} is instrumental in proving that the $q$-difference equations satisfied by $H_{(5,3,0)}(z,q)$, $H_{(4,3,1)}(z,q)$, and $H_{(4,1,3)}(z,q)$ in $S_{11}((a_1,a_2,a_3)|(b_1,b_2,b_3))$ language are elements of $I_{11}$.

Next, we look at the $q$-difference equation satisfied by $H_{(7,1,0)}(z,q)$ from \eqref{eq:Hrec}:
\[ H_{(7, 1, 0)}(z,q) - 
 H_{(7, 0, 1)}(q z,q) - H_{(6, 2, 0)}(q z,q) + (1 - q z) H_{(6, 1, 1)}(q^2 z,q) =0. \]
After the use of \eqref{eq:mod11list} and \eqref{eq:HzTransform}, we see that this $q$-difference equation is equivalent to the following conjectural form
\[S_{11}((0, 1, 1)|(1, 1, 1)) - S_{11}((1, 0, 1)|(1, 1, 1)) - 
 q z S_{11}((2, 1, 1)|(0, 1, 1))=0.\] This is nothing but the relation $R_1^{(1)} (0, 1, 1)|(1, 1, 1)$ of \eqref{R1} given in Lemma~\ref{lemma:recs}. Hence, this relation is also within $I_{11}$ and covered by the relations of $S_{11}((a_1,a_2,a_3)|(b_1,b_2,b_3))$.

As a second explicit example, consider the $q$-difference equation satisfied by $H_{(6,1,1)}(z,q)$, 
\begin{align*}
H_{(6, 1, 1)}(z,q) &- H_{(7, 1, 0)}(q z,q) - H_{(6, 0, 2)}(q z,q)  - H_{(5, 2, 1)}(q z,q)+ (1 - q z) H_{(7, 0, 1)}(q^2 z,q) \\ &  + (1 - q z) H_{(6, 2, 0)}(q^2 z,q) + (1 - q z) H_{(5, 1, 2)}(q^2 z,q)  - (1 - q z) (1 - q^2 z) H_{(6, 1, 1)}(q^3 z,q) = 0.\end{align*}
Using employing \eqref{eq:mod11list} and \eqref{eq:HzTransform}, we see get the (conjecturally) equivalent form \begin{align*}
S_{11}&((0, 1, 1)|(0, 1, 1)) - S_{11}((1, 0, 1)|(0, 1, 1)) - S_{11}((1, 1, 1)|(1, 1, 1))\\ &+ (1- q z)S_{11}((2, 0, 1)|(1, 1, 1))- q z S_{11}((2, 1, 1)|(0, 0, 1)) +  q^2 z(1-q z) S_{11}((3, 1, 1)|(0, 1, 1))=0.
\end{align*}
This relation can be checked to be the side-by-side additions of \[R_1^{(1)}(({{0, 1, 1})|({0, 1, 1}})) - (1 - q z) R_1^{(1)}(({{1, 1, 1})|( {1, 1, 1}})) + 
 q z R_2^{(1)}(({{2, 1, 1})|({-1, 1, 1}})) \in I_{11}.\]
 
 We can one-by-one write down the remaining 8 recurrences, their $S_{11}(\rho|\sigma)$ equivalents, and what combination of \eqref{R1}-\eqref{R40} is equivalent to the functional equations in the $S_{11}(\rho|\sigma)$. This way we prove that these relations are all included in the ideal $I_{11}$. We need to say that these relations gets messier, pages long and not hand-verifiable. Printing these would be a waste of page/paper and instead we keep these in the digital realm for interested readers to check it easily, or print on paper as they wish. To that end, similar to how it was handled in \cite{KR}, we include text files \texttt{M11RecHXYZ\_Explicit.txt} in the ancillary files portion of ArXiv and on the author's website \cite{tinyurl}. Here \texttt{XYZ} is to be replaced by the relevant profile's digits such as \texttt{620} for the profile $(6,2,0)$. One can check that the elements of $I_{11}$ given in these text files are equivalent to the $q$-difference equations \eqref{eq:Hrec} satisfied by $H_{(X,Y,Z)}(z,q)$ after they are translated to $S_{11}(\rho|\sigma)$ form using  \eqref{eq:mod11list}, \eqref{eq:H440exp} and \eqref{eq:HzTransform}. The functional equation names are reflected in the text as \texttt{RX[\{Y\},\{\{a1,a2,a3\},\{b1,b2,b3\}\}]} for \texttt{X} and \texttt{Y} to be replaced by 1 or 2 to denote 
 $R_X^{(Y)}((a_1,a_2,a_3)|(b_1,b_2,b_3))$, or \texttt{RZ[\{\{a1,a2,a3\},\{b1,b2,b3\}\}]} for \texttt{Z} to be replaced by 3 or 4 to denote $R_3((a_1,a_2,a_3)|(b_1,b_2,b_3))$ and $R_4((a_1,a_2,a_3)|(b_1,b_2,b_3))$, respectively. The definitions of these functional equations can be found in Lemma~\ref{lemma:recs} for $m=11$. A guide document that explicitly lists each \texttt{R} functional relation for modulo 10 is given in \texttt{M11R} text file. One also can see that the largest entry within $\rho=(a_1,a_2,a_3)$ and $\sigma=(b_1,b_2,b_3)$ of the relations \eqref{R1}-\eqref{R40} for the modulo 11 case given in the additional documents is $6$. This proves the following theorem and its corollary.
 
 \begin{theorem}\label{thm:n6m11}  Conjecture~\ref{conj:FinIdeal} is correct for $m=11$ and $N=6$.
 \end{theorem}
 
 \begin{corollary}\label{cor:Idealconjm11} Conjecture~\ref{conj:Ideal} is correct for $m=11$.
 \end{corollary}
 
Corollary~\ref{cor:Idealconjm11} is equivalent to the following theorem:
 
 \begin{theorem}\label{thm:m11} The claimed expressions \eqref{eq:mod11list} and \eqref{eq:H440exp} hold.
 \end{theorem}
 
 Observe that Theorem~\ref{thm:m11} adds a new supporting case to Corollary~\ref{KRconjecture}.
 
Now that the main conjectures are proven for the modulus 11 cases, we can specialize $z=1$ and see the 15 sum-product identities coming from the cylindric partitions paradigm.

\begin{theorem}\label{thm:M11sumprod} The following identities hold
\begin{align}\label{eq:M11sumprod}
\sum_{\substack{r_1\geq r_2\geq r_3\geq 0\\s_1\geq s_2\geq s_3\geq0 }} &\frac{q^{r_1^2 - r_1 s_1+ s_1^2+r_2^2 - r_2 s_2+ s_2^2 +r_3^2 + r_3 s_3+ s_3^2}\  p_c(r_1,r_2,r_3,s_1,s_2,s_3,q)}{ (q;q)_{r_1-r_{2}}(q;q)_{r_2-r_{3}}(q;q)_{r_{3}}(q;q)_{s_1-s_{2}}(q;q)_{s_2-s_{3}}(q;q)_{s_{3}}(q;q)_{r_{3}+s_{3}+1}}\\ \nonumber &\hspace{5cm}= \frac{1}{(q;q)_\infty} \frac{1}{\theta(q^{i_1},q^{i_2},q^{i_3},q^{i_4},q^{i_5},q^{i_6},q^{i_7};q^{11})},
\end{align} where the polynomials $p_c(r_1,r_2,r_3,s_1,s_2,s_3,q)$ and the 7-tuples $(i_1,i_2,i_3,i_4,i_5,i_6,i_7)$ for each profile is given in the following table:
\[
\begin{array}{ccc}
\text{Profile  }c   &    p_{c}(r_1,r_2,r_3,s_1,s_2,s_3,q)  & (i_1,i_2,i_3,i_4,i_5,i_6,i_7)    \\ \hline
    (8,0,0)    &    q^{r_1+r_2+r_3+s_1+s_2+s_3}         & (2,3,3,4,4,5,5) \\
    (7,1,0)    &    q^{r_2+r_3+s_1+s_2+s_3}             & (1,2,3,4,4,5,5) \\
    (7,0,1)    &    q^{r_1+r_2+r_3+s_2+s_3}             & (1,2,3,4,4,5,5) \\
    (6,2,0)    &    q^{r_3+s_1+s_2+s_3}                 & (1,2,2,3,4,5,5) \\
    (6,1,1)    &    q^{r_2+r_3+s_2+s_3}(1-q^{r_1+s_1+1})& (1,1,3,3,4,5,5) \\
    (6,0,2)    &    q^{r_1+r_2+r_3+s_3}                 & (1,2,2,3,4,5,5) \\
    (5,3,0)    &    q^{s_1+s_2+s_3}                     & (1,2,2,3,3,4,5) \\
    (5,2,1)    &    q^{r_3+s_2+s_3}(1-q^{r_2+s_1+1})    & (1,1,2,3,4,4,5) \\
    (5,1,2)    &    q^{r_2+r_3+s_3}(1-q^{r_1+s_2+1})    & (1,1,2,3,4,4,5) \\
    (5,0,3)    &    q^{r_1+r_2+r_3}                     & (1,2,2,3,3,4,5) \\
    (4,3,1)    &    q^{s_2+s_3}(1-q^{r_3+s_1+1})        & (1,1,2,3,3,4,5) \\
    (4,2,2)    &    q^{r_3+s_3}(1-q^{r_2+s_2+1})        & (1,1,2,2,4,4,5) \\
    (4,1,3)    &    q^{r_2+r_3}(1-q^{r_1+s_3+1})        & (1,1,2,3,3,4,5) \\
    (3,3,2)    &    q^{s_3}(1-q^{r_3+s_2+1})            & (1,1,2,2,3,5,5) \\ \hline
    (4,4,0)    &    q^{r_1}(q^{s_2+s_3}-q^{r_3+s_1+s_2+s_3+1} + q^{r_1+r_2+r_3+1} )  & (1,2,2,3,3,4,4)
\end{array}
\]
\end{theorem}

In Theorem~\ref{thm:M11sumprod}, we chose to put the profile $(4,4,0)$ related sum-product identity under a line in the table. This is to indicate that this identity is not a direct claim made by combining \eqref{conj:Hconj} with $z=1$ and \eqref{BorodinProd}. We first recovered a formula for $H_{(4,4,0)}(z,q)$ as a combination of $S_{11}((a_1,a_2,a_3)|(b_1,b_2,b_3))$ series and then made this claim. This line also has the added benefit that it aligns us with Kanade--Russell's language as this is the sum-product identity related to the \textit{under-the-line} $H_c(z,q)$ function, which we chose not to directly define.

The sum sides are the expressions \eqref{eq:mod11list} and \eqref{eq:H440exp} with $z=1$ written explicitly using \eqref{eq:Sm1} with $k=4$. The product sides follow from \eqref{eq:HtoF} with $z=1$ followed by \eqref{BorodinProd}. The product related to the first profile, $(8,0,0)$, on the table is presented in the introduction as Theorem~\ref{th:mod11ex}.

Observe that the products that appear on the right-hand side of \eqref{eq:M11sumprod} related to the profiles $(c_1,c_2,c_3)$ and  $(c_1,c_3,c_2)$ are the same. The symmetry for the generating functions have been observed and noted before, for example in \cite[Corollary~2.2]{CDU}. This symmetry is visible on the sum side of Theorem~\ref{thm:M11sumprod} too. One can get the ``other" sum by merely replacing the variable `$r$'s and `$s$'s. Note that this is a byproduct of setting $z=1$ and this similarity does not exist on the sum side for generic $z$. In that light, this theorem consisting of 15 sum-product identities actually provide a total of 10 essentially unique sum-product identities.

We also note that among these identities the ones related to profiles $(8,0,0)$, $(6,1,1)$, $(4,2,2)$, $(3,3,2)$ are the $i=1,\dots,4$ cases of (5.28), respectively, and $(5,3,0)$ and $(6,2,0)$ are the $\sigma=0$ and $1$ cases of (5.29), respectively, of \cite[Theorem~5.3]{ASW99}.

\section{Modulo 13 Identities}\label{sec:M13}

Similar to Section~\ref{sec:M11}, we start by listing the explicit claims of Conjecture~\ref{conj:Hconj} for the modulus  $m=13$ family.

\begin{equation}\label{eq:mod13list}
\begin{array}{ll}
H_{(10, 0, 0)}(z,q) &= S_{13}((1, 1, 1)|(1, 1, 1)),\\
H_{(9, 1, 0)}(z,q) &= S_{13}((0, 1, 1)|(1, 1, 1)),\\
H_{(9, 0, 1)}(z,q) &= 
 S_{13}((1, 1, 1)|(0, 1, 1)) - q (1 - z) S_{13}((2, 1, 1)|(1, 1, 1)),\\
H_{(8, 2, 0)}(z,q) &= S_{13}((0, 0, 1)|(1, 1, 1)),\\
H_{(8, 1, 1)}(z,q) &= 
 S_{13}((0, 1, 1)|(0, 1, 1)) - q S_{13}((1, 1, 1)|(1, 1, 1)),\\
H_{(8, 0, 2)}(z,q) &= 
 S_{13}((1, 1, 1)|(0, 0, 1)) - q (1 - z) S_{13}((2, 1, 1)|(0, 1, 1)),\\ 
H_{(7, 3, 0)}(z,q) &= S_{13}((0, 0, 0)|(1, 1, 1)),\\
H_{(7, 2, 1)}(z,q) &= 
 S_{13}((0, 0, 1)|(0, 1, 1)) - q S_{13}((0, 1, 1)|(1, 1, 1)),\\ 
H_{(7, 1, 2)}(z,q) &= 
 S_{13}((0, 1, 1)|(0, 0, 1)) - q S_{13}((1, 1, 1)|(0, 1, 1)),\\ 
H_{(7, 0, 3)}(z,q) &= 
 S_{13}((1, 1, 1)|(0, 0, 0)) - q (1 - z) S_{13}((2, 1, 1)|(0, 0, 1)),\\ 
H_{(6, 3, 1)}(z,q) &= 
 S_{13}((0, 0, 0)|(0, 1, 1)) - q S_{13}((0, 0, 1)|(1, 1, 1)),\\ 
H_{(6, 2, 2)}(z,q) &= 
 S_{13}((0, 0, 1)|(0, 0, 1)) - q S_{13}((0, 1, 1)|(0, 1, 1)),\\ 
H_{(6, 1, 3)}(z,q) &= 
 S_{13}((0, 1, 1)|(0, 0, 0)) - q S_{13}((1, 1, 1)|(0, 0, 1)),\\ 
H_{(5, 3, 2)}(z,q) &= 
 S_{13}((0, 0, 0)|(0, 0, 1)) - q S_{13}((0, 0, 1)|(0, 1, 1)),\\ 
H_{(5, 2, 3)}(z,q) &= 
 S_{13}((0, 0, 1)|(0, 0, 0)) - q S_{13}((0, 1, 1)|(0, 0, 1)),\\ 
H_{(4, 3, 3)}(z,q) &= 
 S_{13}((0, 0, 0)|(0, 0, 0)) - q S_{13}((0, 0, 1)|(0, 0, 1)).
\end{array}
\end{equation}

There are six profiles that are not covered by Conjecture~\ref{conj:Hconj}. Once again, using \eqref{eq:Hrec} explicit claims for the normalized generating functions related to the number of cylindric partitions with these profiles can be recovered. We make the claims in the following succession.

First we look at the $q$-difference equation \eqref{eq:Hrec} that $H_{(7,3,0)}$: \begin{equation}\label{eq:recH730} H_{(7, 3, 0)}(z,q)- 
 H_{(7, 2, 1)}(q z,q) - H_{(6, 4, 0)}(q z,q)  + (1 - q z) H_{(6, 3, 1)}(q^2 z,q)= 0.\end{equation} By writing the $S_{13}((a_1,a_2,a_3)|(b_1,b_2,b_3))$ equivalents for the functions in \eqref{eq:mod13list} and using \eqref{eq:HzTransform}, we get

\begin{align}\label{eq:defH640}
   H_{(6, 4, 0)}(z,q) &= S_{13}((-1, 0, 0)|(1, 1, 1)) - S_{13}((0, 0, 1)|(0, 1, 1)) + 
 q S_{13}((0, 1, 1)|(1, 1, 1))\\\nonumber &+ (1 - z) S_{13}((1, 0, 0)|(0, 1, 1)) - q (1 - z)S_{13}((1, 0, 1)|(1, 1, 1)).
\end{align}

Note that we did not use the $q$-difference equation of $H_{(6,4,0)}(z,q)$ to make a claim for its formula. In Section~\ref{sec:M11}, there was only a single missing formula. That allowed us to use the $q$-difference equation for that very function and get a formula in $S_{11}((a_1,a_2,a_3)|(b_1,b_2,b_3))$'s with no backwards shifts (i.e. $z\mapsto z/q$, which also reflects as negative indices in the first variable $a_1$). This may not be possible in general. The $q$-difference equation $H_{(6, 4, 0)}(z,q)$ satisfies is \begin{equation}\label{eq:recH640} H_{(6, 4, 0)}(z,q)-  H_{(6, 3, 1)}(q z,q)  - H_{(5, 5, 0)}(q z,q) +(1 - q z) H_{(5, 4, 1)}(q^2 z,q) = 0.\end{equation} The conjectural formulas \eqref{eq:mod13list} does not cover $ H_{(5, 5, 0)}(z,q)$. Hence, we cannot fully translate $H_{(6, 4, 0)}(z,q)$ to a formula made up of $S_{13}((a_1,a_2,a_3)|(b_1,b_2,b_3))$ series. Nevertheless, as also noted in \cite{KR}, we can recover formulas for all the missing functions using other recurrences and backwards shifts in $a_1$. 

In fact, the recurrence \eqref{eq:recH640} and \eqref{eq:defH640} can be put together to claim a formula for $H_{(5,5,0)}(z,q)$. After the similar considerations we claim \begin{align}\label{eq:defH550}
H_{(5,5,0)}(z,q) &= S_{13}((-2, 0, 0)|(1, 1, 1)) - S_{13}((-1, 0, 1)|(0, 1, 1)) + 
 q S_{13}((-1, 1, 1)|(1, 1, 1))\\ \nonumber&+ (1-z/q-z) S_{13}((0, 0, 0)|(0, 1, 1)) - (q - z - 
    q z) S_{13}((0, 0, 1)|(1, 1, 1)) \\\nonumber &- 
 q (1 - z) z S_{13}((1, 0, 0)|(1, 1, 1)) - (1 - 
    z) S_{13}((1, 0, 1)|(0, 0, 1)) \\ \nonumber &+q (1 - z) S_{13}((1, 1, 1)|(0, 1, 1)) + (1 - z) (1 - 
    q z) S_{13}((2, 0, 0)|(0, 0, 1)) \\ \nonumber &-
 q (1 - z) (1 - q z) S_{13}((2, 0, 1)|(0, 1, 1)).
\end{align}
We point out that the coefficients of the claimed $H_{(5,5,0)}(z,q)$ formula now can be seen to have a Laurent polynomial. This is a byproduct of the backwards shifts in $z$.

Using the $q$-difference equation for $H_{(6,3,1)}(z,q)$,
\begin{align}\label{eq:recH631}
H_{(6, 3, 1)}(z,q)&- H_{(7, 3, 0)}(q z,q)- H_{(6, 2, 2)}(q z,q)- H_{(5, 4, 1)}(q z,q) + (1 - q z) H_{(7, 2, 1)}( q^2 z,q)\\ \nonumber &+ (1 - q z) H_{(6, 4, 0)}(q^2 z,q) + (1 - q z) H_{(5, 3, 2)}(q^2 z,q)   - (1 - q z) (1 - q^2 z) H_{(6, 3, 1)}( q^3 z,q) = 0,
\end{align}
   \eqref{eq:mod13list} and \eqref{eq:defH640} we claim that 
\begin{align} \label{eq:defH541}
 H_{(5, 4, 1)}(z,q) &=  S_{13}((-1, 0, 0)|(0, 1, 1)) - q S_{13}((-1, 0, 1)|(1, 1, 1)) - z S_{13}((0, 0, 0)|(1, 1, 1))\\
 \nonumber &- S_{13}((0, 0, 1)|(0, 0, 1))  + q S_{13}((0, 1, 1)|(0, 1, 1)) +(1-z) S_{13}((1, 0, 0)|(0, 0, 1))\\ 
 \nonumber & - q(1-z) S_{13}((1, 0, 1)|(0, 1, 1)).
\end{align}

Using the $q$-difference equation for $H_{(5,3,2)}(z,q)$,
\begin{align}\label{eq:recH532}
H_{(5, 3, 2)}(z,q) &- H_{(6, 3, 1)}(q z,q)-  H_{(5, 2, 3)}(q z,q)- H_{(4, 4, 2)}(q z,q)+ (1 - q z) H_{(6, 2, 2)}(q^2 z,q)\\ \nonumber & + (1 - q z) H_{(5, 4, 1)}(q^2 z,q)  + (1 - q z) H_{(4, 3, 3)}(q^2 z,q)- (1 - q z) (1 - q^2 z) H_{(5, 3, 2)}(q^3 z,q) = 0,
\end{align}
   \eqref{eq:mod13list} and \eqref{eq:defH541} we claim that 
\begin{align} \label{eq:defH442}
 H_{(4, 4, 2)}(z,q) &=  S_{13}((-1, 0, 0)|(0, 0, 1)) - q S_{13}((-1, 0, 1)|(0, 1, 1)) -z S_{13}((0, 0, 0)|(0, 1, 1))\\ \nonumber & - S_{13}((0, 0, 1)|(0, 0, 0)) + q z S_{13}((0, 0, 1)|(1, 1, 1)) + q S_{13}((0, 1, 1)|(0, 0, 1))\\ \nonumber & + (1-z)S_{13}((1, 0, 0)|(0, 0, 0)) -  q z(1-z) S_{13}((1, 0, 0)|(1, 1, 1)) \\ 
 \nonumber & - q(1-z) S_{13}((1, 0, 1)|(0, 0, 1)).
\end{align}

Then, by the $q$-difference equation for $H_{(5,2,3)}(z,q)$,
\begin{align}\label{eq:recH523}H_{(5, 2, 3)}(z,q)&- H_{(6, 2, 2)}(q z,q)- H_{(5, 1, 4)}(q z,q)-H_{(4, 3, 3)}(q z,q)+ (1 - q z) H_{(6, 1, 3)}(q^2 z,q)\\ \nonumber &+ (1 - q z) H_{(5, 3, 2)}(q^2 z,q)+ (1 - q z) H_{(4, 4, 2)}(q^2 z,q)   - (1 - q z) (1 - q^2 z) H_{(5, 2, 3)}(q^3 z,q)    = 0,\end{align}
   together with \eqref{eq:mod13list} and \eqref{eq:defH442} we claim that 
\begin{align} \label{eq:defH514}
 H_{(5,1,4)}(z,q) &=  S_{13}((-1, 0, 1)|(0, 0, 0)) - q S_{13}((-1, 1, 1)|(0, 0, 1)) - S_{13}((0, 0, 0)|(0, 0, 0)) \\ 
 \nonumber &+ (1 - z) S_{13}((0, 0, 0)|(0, 0, 1)) - (1 - q) S_{13}((0, 0, 1)|(0, 0, 1)) \\ 
 \nonumber &- q (1 - z) S_{13}((0, 0, 1)|(0, 1, 1)) +  q S_{13}((0, 1, 1)|(0, 1, 1)) \\ 
 \nonumber &+ (1 - z) S_{13}((1, 0, 0)|(0, 0, 1)) - q (1 - z) z S_{13}((1, 0, 0)|(0, 1, 1)) \\ 
 \nonumber &- (1 - z) S_{13}((1, 0, 1)|(0, 0, 0)) - q (1 - z) S_{13}((1, 0, 1)|(0, 1, 1)) \\ 
 \nonumber &+ q^2 (1 - z) z S_{13}((1, 0, 1)|(1, 1, 1)) + (1 - z) S_{13}((1, 1, 1)|(0, 0, 0)) \\ 
 \nonumber &+ q (1 - z) S_{13}((1, 1, 1)|(0, 0, 1)) + (1 - z) (1 - q z) S_{13}((2, 0, 0)|(0, 0, 0)) \\ 
 \nonumber &- q^2 z (1 - z)(1 - q z) S_{13}((2, 0, 0)|(1, 1, 1)) - (1 - z) (1 - q z) S_{13}((2, 0, 1)|(0, 0, 0))\\
 \nonumber &- q (1 - z) (1 - q z) S_{13}((2, 0, 1)|(0, 0, 1)) - q^2 z (1 - z)  S_{13}((2, 1, 1)|(0, 0, 1)).
\end{align}

Finally, by replacing the formulas in \eqref{eq:mod13list} and \eqref{eq:defH514} in \begin{equation}\label{eq:recH604}H_{(6, 0, 4)}(z,q)-  H_{(7, 0, 3)}(q z,q) -H_{(5, 1, 4)}(q z,q)  + (1 - q z) H_{(6, 1, 3)}(q^2 z,q) =0\end{equation} we conjecture that \begin{align}
\label{eq:defH604}
 H_{(6,0,4)}(z,q) &= S_{13}((0, 0, 1)|(0, 0, 0)) - q S_{13}((0, 1, 1)|(0, 0, 1)) - S_{13}((1, 0, 0)|(0, 0, 0)) \\
 \nonumber &+ (1 - q z) S_{13}((1, 0, 0)|(0, 0, 1)) - (1 -  q) S_{13}((1, 0, 1)|(0, 0, 1)) \\ 
 \nonumber &-  q (1 - q z) S_{13}((1, 0, 1)|(0, 1, 1)) +  q S_{13}((1, 1, 1)|(0, 1, 1)) \\ 
 \nonumber &+ (1 - q z) S_{13}((2, 0, 0)|(0, 0, 1)) -  q^2 z (1 - q z) S_{13}((2, 0, 0)|(0, 1, 1)) \\ 
 \nonumber &- (1 - q z) S_{13}((2, 0, 1)|(0, 0, 0)) -  q (1 - q z) S_{13}((2, 0, 1)|(0, 1, 1)) \\ 
 \nonumber &+  q^3 z (1 - q z) S_{13}((2, 0, 1)|(1, 1, 1)) +  S_{13}((2, 1, 1)|(0, 0, 0)) \\ 
 \nonumber &+  q (1 - q z) S_{13}((2, 1, 1)|(0, 0, 1)) + (1 - q z) (1 - q^2 z) S_{13}((3, 0, 0)|(0, 0, 0)) \\ 
 \nonumber &-  q^3 z (1 - q z) (1 - q^2 z) S_{13}((3, 0, 0)|(1, 1, 1)) - (1 - q z) (1 - q^2 z) S_{13}((3, 0, 1)|(0, 0, 0)) \\ 
 \nonumber &- q (1 - q z) (1 - q^2 z) S_{13}((3, 0, 1)|(0, 0, 1)) - q^3 z (1 - q z) S_{13}((3, 1, 1)|(0, 0, 1)).
  \end{align}
  
We can prove that the later claimed $H_{(6,4,0)}(z,q)$, $H_{(5,5,0)}(z,q)$, $H_{(5,4,1)}(z,q)$, $H_{(5,3,2)}(z,q)$, $H_{(5,2,3)}(z,q)$, and $H_{(6,0,4)}(z,q)$  the initial conditions $H_c(0,q)=1/(q;q)_\infty$ and $H_c(z,0)=1$ in the succession from \eqref{eq:recH730}, \eqref{eq:recH640}, \eqref{eq:recH631}, \eqref{eq:recH532}, \eqref{eq:recH523}, and \eqref{eq:recH604}, respectively. To prove  the  $H_c(0,q)=1/(q;q)_\infty$ initial condition we need to first shift $z\mapsto z/q$ in all but the last of the functional equations.

The $q$-difference equations for $H_{(10,0,0)}(z,q)$, $H_{(9,0,1)}(z,q)$, $H_{(8,0,2)}(z,q)$, and $H_{(7,0,3)}(z,q)$ becomes tautologies once translated into $S_{13}$ form using \eqref{eq:mod13list} and \eqref{eq:HzTransform}. The $q$-difference equations for $H_{(7,3,0)}(z,q)$, $H_{(6,4,0)}(z,q)$, $H_{(6,3,1)}(z,q)$, $H_{(5,3,2)}(z,q)$, $H_{(5,2,3)}(z,q)$, and $H_{(6,0,4)}(z,q)$ are the recurrences used to define the missing $H_c(z,q)$ functions in the modulo 13 family (see \eqref{eq:recH730}, \eqref{eq:recH640}, \eqref{eq:recH631}, \eqref{eq:recH532}, \eqref{eq:recH523}, and \eqref{eq:recH604}, resp.). Hence, these equations also trivializes once the relevant functions are written in their claimed $S_{13}$ forms using \eqref{eq:mod13list}, \eqref{eq:defH640}, \eqref{eq:defH550}, \eqref{eq:defH541}, \eqref{eq:defH442}, \eqref{eq:defH514}, and \eqref{eq:defH604} together with \eqref{eq:HzTransform}. 

After the considerations above, we end up with 10 non-trivial coupled $q$-difference equations to prove. Showing that the $q$-difference equations' in the claimed $S_{13}((a_1,a_2,a_3)|(b_1,b_2,b_3))$ belong to the ideal $I_{13}$, which is generated by the relations of $S_{13}((a_1,a_2,a_3)|(b_1,b_2,b_3))$s (see Lemma~\ref{lemma:recs}), is done by the method outlined in Section~\ref{sec:proofmethod}.  Explicit linear combination of \eqref{R1}-\eqref{R40} equivalents of these 12 functional equations in $S_{13}$ form can, once again, be found in the ancillary files portion of ArXiv and on the author's website \cite{tinyurl} under the file names \texttt{M13RecHXYZ\_Explicit.txt}. Here \texttt{XYZ} is to be replaced by the relevant profile's digits such as \texttt{910} for the profile $(9,1,0)$. One can check that the elements of $I_{13}$ given in these text files are equivalent to the $q$-difference equations \eqref{eq:Hrec} satisfied by $H_{(X,Y,Z)}(z,q)$ after they are translated to $S_{11}(\rho|\sigma)$ form using \eqref{eq:mod13list}, \eqref{eq:defH640}, \eqref{eq:defH550}, \eqref{eq:defH541}, \eqref{eq:defH442}, \eqref{eq:defH514}, and \eqref{eq:defH604} and \eqref{eq:HzTransform}. The recurrence names are reflected in the text as \texttt{RX[\{Y\},\{\{a1,a2,a3\},\{b1,b2,b3\}\}]} for  \texttt{X} and \texttt{Y} to be replaced by 1 or 2 to denote  $R_X^{(Y)}((a_1,a_2,a_3)|(b_1,b_2,b_3))$, or \texttt{RZ[\{\{a1,a2,a3\},\{b1,b2,b3\}\}]} for \texttt{Z} to be replaced by 3 or 4 to denote $R_3((a_1,a_2,a_3)|(b_1,b_2,b_3))$ and $R_4((a_1,a_2,a_3)|(b_1,b_2,b_3))$, respectively. Finally, a guide document that explicitly lists each \texttt{R} functional relation for modulo 10 is given in \texttt{M13R} text file. 

This tedious, error prone and impossible by hand calculation proves the following theorem and its corollary.
 
 \begin{theorem}\label{thm:n6m13}  Conjecture~\ref{conj:FinIdeal} is correct for $m=13$ and $N=6$.
 \end{theorem}
 
 \begin{corollary}\label{cor:Idealconjm13} Conjecture~\ref{conj:Ideal} is correct for $m=13$.
 \end{corollary}
 
Corollary~\ref{cor:Idealconjm13} is equivalent to the following theorem:
 
 \begin{theorem}\label{thm:m13} The claimed expressions of \eqref{eq:mod13list},  \eqref{eq:defH640}, \eqref{eq:defH550}, \eqref{eq:defH541}, \eqref{eq:defH442}, \eqref{eq:defH514}, and \eqref{eq:defH604} hold.
 \end{theorem}
 
 As before, Theorem~\ref{thm:m13} adds another new witness to Corollary~\ref{KRconjecture}, and increases our confidence in it.
 
Now that the main conjectures are proven for the modulus 13 cases, we can set $z=1$ and see the 22 sum-product identities coming from the cylindric partitions paradigm.

\begin{theorem}\label{thm:M13sumprod} The following identities hold
\begin{align}\label{eq:M13sumprod}
\sum_{\substack{r_1\geq r_2\geq r_3\geq 0\\s_1\geq s_2\geq s_3\geq0 }} &\frac{q^{r_1^2 - r_1 s_1+ s_1^2+r_2^2 - r_2 s_2+ s_2^2 +r_3^2 - r_3 s_3+ s_3^2}\  p_c(r_1,r_2,r_3,s_1,s_2,s_3,q)}{ (q;q)_{r_1-r_{2}}(q;q)_{r_2-r_{3}}(q;q)_{r_{3}}(q;q)_{s_1-s_{2}}(q;q)_{s_2-s_{3}}(q;q)_{s_{3}}(q;q)_{r_{3}+s_{3}+1}}\\ \nonumber &\hspace{5cm}= \frac{1}{(q;q)_\infty} \frac{1}{\theta(q^{i_1},q^{i_2},q^{i_3},q^{i_4},q^{i_5},q^{i_6},q^{i_7},q^{i_8},q^{i_9};q^{13})},
\end{align} where the polynomials $p_c(r_1,r_2,r_3,s_1,s_2,s_3,q)$ and the 9-tuples $(i_1,i_2,i_3,i_4,i_5,i_6,i_7,i_8,i_9)$ for each profile is given in the following table:
\[
\begin{array}{ccc}
\text{Profile  }c   &    p_{c}(r_1,r_2,r_3,s_1,s_2,s_3,q)  & (i_1,i_2,i_3,i_4,i_5,i_6,i_7,i_8,i_9)    \\ \hline
    (10,0,0)        &    q^{r_1+r_2+r_3+s_1+s_2+s_3}         & (2,3,3,4,4,5,5,6,6) \\
    (9,1,0)         &    q^{r_2+r_3+s_1+s_2+s_3}             & (1,2,3,4,4,5,5,6,6) \\
    (9,0,1)         &    q^{r_1+r_2+r_3+s_2+s_3}             & (1,2,3,4,4,5,5,6,6) \\
    (8,2,0)         &    q^{r_3+s_1+s_2+s_3}                 & (1,2,2,3,4,5,5,6,6) \\
    (8,1,1)         &    q^{r_2+r_3+s_2+s_3}(1-q^{r_1+s_1+1})& (1,1,3,3,4,5,5,6,6) \\
    (8,0,2)         &    q^{r_1+r_2+r_3+s_3}                 & (1,2,2,3,4,5,5,6,6) \\
    (7,3,0)         &    q^{s_1+s_2+s_3}                     & (1,2,2,3,3,4,5,6,6) \\
    (7,2,1)         &    q^{r_3+s_2+s_3}(1-q^{r_2+s_1+1})    & (1,1,2,3,4,4,5,6,6) \\
    (7,1,2)         &    q^{r_2+r_3+s_3}(1-q^{r_1+s_2+1})    & (1,1,2,3,4,4,5,6,6) \\
    (7,0,3)         &    q^{r_1+r_2+r_3}                     & (1,2,2,3,3,4,5,6,6) \\
    (6,3,1)         &    q^{s_2+s_3}(1-q^{r_3+s_1+1})        & (1,1,2,3,3,4,5,5,6) \\
    (6,2,2)         &    q^{r_3+s_3}(1-q^{r_2+s_2+1})        & (1,1,2,2,4,4,5,5,6) \\
    (6,1,3)         &    q^{r_2+r_3}(1-q^{r_1+s_3+1})        & (1,1,2,3,3,4,5,5,6) \\
    (5,3,2)         &    q^{s_3}(1-q^{r_3+s_2+1})            & (1,1,2,2,3,4,5,5,6) \\
    (5,2,3)         &    q^{s_3}(1-q^{r_3+s_2+1})            & (1,1,2,2,3,4,5,5,6) \\
    (4,3,3)         &    (1-q^{r_3+s_3+1})                   & (1,1,2,2,3,3,5,6,6) \\ \hline
    (6,4,0)         &    q^{s_2+s_3}(q^{-r_1 + s_1} - q^{r_3}+ q^{r_2 + r_3 + s_1 + 1}) & (1,2,2,3,3,4,4,5,6)   \\[-1.5ex]\\
    (6,0,4)         &   
    \begin{array}{l}
q^{r_3} - q^{r_2 + r_3 + s_3 + 1} - q^{r_1}+ q^{2 r_1 + r_2 + r_3}+ q^{r_1 + r_2 + r_3 + s_2 + s_3 + 1}\\
+ (1 - q) q^{r_1+s_3} (1-q^{r_3} - q^{r_3+s_3+1}) \\- (1-q)q^{2r_1}(q^{r_3} - q^{s_3} 
- q^{ r_2 + r_3 + s_3+1} + q^{ r_1 + r_2 + r_3 + s_3+3}\\\hspace{2cm}+q^{ s_2 + s_3+2} + q^{ r_3 + s_2 + s_3+1} - q^{ r_3 + s_1 + s_2 + s_3+3})\\+(1-q)(1-q^2)q^{r_3}(1-q^{r_3} - q^{r_3+s_3+1}+q^{s_1+s_2+s_3+3})\\
    \end{array}     & (1,2,2,3,3,4,4,5,6)  \\[-1.5ex]\\
    (5,5,0)         &    \begin{array}{l}q^{-2 r_1 + s_1 + s_2 + s_3} - q^{-r_1 + r_3 + s_2 + s_3}- q^{s_2 + s_3 - 1} +  q^{r_3 + s_1 + s_2 + s_3} \\+  q^{-r_1 + r_2 + r_3 + s_1 + s_2 + s_3 + 1}  \end{array}                           &      (1,2,2,3,3,4,4,5,5)                 \\[-1.5ex]\\
    (5,4,1)         &         \begin{array}{l}q^{-r_1 + s_2 + s_3} - q^{-r_1 + r_3 + s_1 + s_2 + s_3 + 1} - 
 q^{s_1 + s_2 + s_3}- q^{r_3 + s_3}\\  + q^{r_2 + r_3 + s_2 + s_3 + 1}
    \end{array}&      (1,1,2,3,3,4,4,5,6) \\[-1.5ex]\\
    (5,1,4)         &  q^{-r_1 + r_3} - q^{-r_1 + r_2 + r_3 + s_3 + 1} + 
 q^{r_2 + r_3 + s_2 + s_3 + 1}- (1 - q) q^{r_3 + s_3} -1       & (1,1,2,3,3,4,4,5,6)\\[-1.5ex]\\
    (4,4,2)         &    \begin{array}{l}q^{-r_1 + s_3} - q^{-r_1 + r_3 + s_2 + s_3 + 1} - q^{s_2 + s_3} - q^{r_3}+  q^{r_3 + s_1 + s_2 + s_3 + 1}\\  + q^{r_2 + r_3 + s_3 + 1}
    \end{array}&      (1,1,2,2,3,4,4,6,6)\\[-1.5ex]\\
\end{array}
\]
\end{theorem}

Once we ignore the symmetries between variables $r$ and $s$, Theorem~\ref{thm:M13sumprod} proves 16 essentially unique sum-product identities. It can easily be seen that within the under-the-line identities, we do not see these symmetries. The product related to the first profile, $(10,0,0)$, on the table is presented in the introduction as Theorem~\ref{th:mod13ex}.

We also note that among these identities the ones related to profiles $(10,0,0)$, $(8,1,1)$, $(6,2,2)$, $(4,3,3)$ are the $i=1,\dots,4$ cases of (5.22), respectively, and $(7,3,0)$ and $(8,2,0)$ are the $\sigma=0$ and $1$ cases of (5.23), respectively, of \cite[Theorem~5.1]{ASW99}.

\section{Future Directions}\label{sec:future}

There are many mathematical questions that arose from the recent studies on cylindric partitions. It is relevant to mention some of the future directions we plan to pursue. 

The approach outlined in \cite{KR} and in this paper attempts to prove sum-representations for all the normalized generating function $H_c(z,q)$ in one stroke for any fixed $|c|$ where $\#(c)=3$. The proof requires hefty calculations after the under-the-line sums are recovered. Then by setting $z=1$ and using \eqref{BorodinProd}, we prove sum-product identities for all profiles within a cylindric partition system for a fixed modulus, again in one stroke. Therefore, to prove $A_2$ Rogers--Ramanujan identities we first prove a more general and more complicated combinatorial connection with a free variable $z$. The success of this method depends on the completion of these calculations, which is virtually impossible by hand. 

Warnaar \cite{WarPriv} mentioned that he build the necessary theory of the Bailey machinery for profiles with 3 parts. This machinery will allow us to prove one sum-product identity at a time. This is wonderful to hear and a great advancement in mathematics. Sadly, it comes with its own short-comings. Warnaar acknowledged that this Bailey machinery can not prove any under-the-line identity at the moment. It can only find the sum-product relation related to the $z=1$ specializations of Conjecture~\ref{eq:Hconj}. This is similar to the situation of the original Andrews--Schilling--Warnaar paper, where for example at the modulo 7 case the Bailey machinery there couldn't reach the under-the-line identity related to the profile $(2,2,0)$, which was later proven in \cite{CW19}.

Be that as it may, we plan to investigate ways to simplify calculations necessary to prove the identities as a whole in one stroke for the free $z$ case by adding the extra information we gather from Warnaar's results. At the very least, for the $z=1$ specialization, we should pursue ways to prove under-the-line identities using the Bailey-machinery-proven over-the-line identities.

There are other sum-product identities that are not visible through the cylindric partitions paradigm.These identities do not have a related cylindric partition profiles attached to them either. Similar to the under-the-line identities, we discover and prove these sum representations using the proven relations in the cylindric partitions system. For example, there are the following two modulo 10 examples similar to \eqref{eq:ASW500}:
\begin{align}
\label{A1}
    \sum_{\substack{r_1\geq r_2\geq0\\s_1\geq s_2\geq0} }\frac{q^{r_1^2 - r_1 s_1+ s_1^2+r_2^2 - r_2 s_2+ s_2^2}\, q^{s_1+s_2}(1+q^{r_1+r_2+1})}{ (q;q)_{r_1-r_{2}}(q;q)_{s_1-s_{2}}(q;q)_{r_{2}}(q;q)_{s_{2}}(q;q)_{r_{2}+s_{2}+1}}&= \frac{1}{(q;q)_\infty} \frac{1}{\theta(q,q,q^{3},q^{4},q^{4},q^{4};q^{10})},\\
\label{A2}    \sum_{\substack{r_1\geq r_2\geq0\\s_1\geq s_2\geq0} }\frac{q^{r_1^2 - r_1 s_1+ s_1^2+r_2^2 - r_2 s_2+ s_2^2}\, q^{s_1+s_2}(1-q^{r_1+r_2+1})}{ (q;q)_{r_1-r_{2}}(q;q)_{s_1-s_{2}}(q;q)_{r_{2}}(q;q)_{s_{2}}(q;q)_{r_{2}+s_{2}+1}}&= \frac{1}{(q;q)_\infty} \frac{1}{\theta(q^{2},q^{2},q^{2},q^{3},q^{3},q^{3};q^{10})}.
\end{align}
All the products associated to principal characters of modulo 10 $A_2$ Rogers--Ramanujan identities are covered by the products that appear in \eqref{BorodinProd}. The identities \eqref{A1} and \eqref{A2} are outside of this system and appear, so to say, on the dark-side of the cylinder. We hope to find a cylindric partition interpretation of these identities in the future. Nevertheless, we plan to present the proofs of these theorems using $q$-theoretic means in an upcoming paper. 

It is still highly relevant to find manifestly positive sum representations for any one of the identities mentioned here. We are looking for ways to see the positivity of the series coefficients. In \cite{ASW99}, Andrews--Schilling--Warnaar suggests applying hypergeometric transformations to eliminate the $(q;q)_\infty$ factor that appear in the identities (such as \eqref{eq:ASW500}) to get a manifestly positive representation. That suggestion is limited and might not be widely applicable, especially for the under-the-line identities. 

In the study of symmetric cylindric partitions \cite{BU} another two fundamental modulo 8 partition theoretic identity families, namely G\"ollnitz--Gordon and little G\"ollnitz identities, showed up. The G\"ollnitz--Gordon identities are known to be related to the level 2 modules of affine Lie algebra $A_5^{(2)}$ \cite{KGG}. This raises new questions of whether,  similar to the symmetric partitions paradigm, we can also relate symmetric cylindric partitions to character formulas of some affine Lie algebras. The product formula analogous to \eqref{BorodinProd} for the count of symmetric cylindric partitions' is present in \cite{BU}. At the moment, the relation of these products' to affine Lie algebra character formulas are fuzzy, and there are no general conjectural series representations for symmetric cylindric partitions either. We plan to study these objects further.

Finally, we plan to pursue sum representations of any generating functions for cylindric partitions with profiles of more than 3 parts. The product representation \eqref{BorodinProd} and the functional equations \eqref{CorteelRec} apply regardless of the size and length of the profiles. So far, we are only able to prove and conjecture sum representations for the profiles with up to 3 parts. 

\section{Comments on Computations}\label{sec:computations}

In the computerized proofs of \cite{CDU}, we make extensive use of \cite{qFunctions} and \cite{HolonomicFunctions}. Those proofs had three main steps. Finding a recurrence relation (over the exponent of $z$) for claimed sum formulas of the (normalized) generating functions of cylindric partitions, uncoupling the $q$-difference equation system laid out by the \eqref{CorteelRec} to get a recurrence satisfied by the coefficient of the $z$'s in the true generating functions of cylindric partitions, comparing recurrences (taking greatest common divisors of recurrences as operators if needed) and showing that both sequences satisfy the same recurrences with the same initial conditions. Once the critical mass of proved identities were reached the rest of the identities were shown by series manipulations guided by \eqref{CorteelRec}. That way we showed that all the claimed sum and the true combinatorial generating function were the same. This proof required two hefty algorithms, namely Creative Telescoping algorithm and Gr\"obner bases calculations, to find the recurrence of a given hypergeometric sum dependent of a discrete variable and to uncouple a coupled system of recurrences, respectively. 

We tried using the same method to prove some claims Warnaar \cite{War21} made for cylindric partitions with 3 part profiles where the modulus is not divisible by 3. Then we quickly saw that the Creative Telescoping calculations were not terminating (in any definition of reasonable time). This is due to the increasing number of nested summations in these conjectures. However, uncoupling of recurrences could still be performed. 

Kanade--Russell's approach \cite{KR} to prove that the claimed series representations for the bivariate generating functions of cylindric partitions are the true generating functions is a fresh take on things. It is somehow backwards compared to the proofs of \cite{CDU}, in the sense that we first extend our conjectural identities using the explicit conjectures of Conjecture~\ref{eq:Hconj} and series manipulations, then prove all these conjectural identities by showing that the coupled relations are satisfied and that we still satisfy the initial conditions. The key idea of reducing coupled $q$-difference equation with the functional relations of the claimed hypergeometric sums was also used in \cite{Chern} in a different context. Moreover, this approach replaces (the old bottle-neck) Creative Telescoping with the contiguous relations of Lemma~\ref{lemma:recs}. However, rewriting the coupled relations of \eqref{CorteelRec} in the new language as a linear combination of terms in the ideal $I_m$ (see Section~\ref{sec:proofmethod}) with coefficients in $\mathbb{Z}(\!(q,z)\!)$ is highly non-trivial. Kanade \cite{Kpriv} mentioned that they found these linear combinations by first making an ansatz for a single case at a time and then solving for undetermined coefficients. The identification of the minimal necessary ansatz is impossible. They also mentioned that each hard-case proof of modulo 10 calculations took about 8 hours to terminate on a home computer. With the matrix reduction approach of this paper, we are order of 2 faster in the modulo 10 cases. This is basically because once we reduce a matrix, we can use it repeadetly for all the functional relations, whereas the previous approach needs to make a single ansatz and solve if for all cases individually. It is with this speed upgrade that we could prove the new modulo 11 and modulo 13 cases. On the other hand, modulo 9 and modulo 12 cases are still open. This is likely due to the extra degree of complication the $q$-binomial coefficients in \eqref{eq:S0}'s introduce. As the order of the recurrences the $S_m(\rho|\sigma)$ satisfy increases, the systems we need to reduce also become larger.

 Mathematica's Gaussian elimination function \texttt{RowReduce} is adamant in calculating the reduced row echelon form of matrices. This is not only not necessary, it also overcomplicates the calculations by introducing large rational function expressions for upper triangular coefficients. This forced us to implementing our own Gaussian elimination algorithm within the Mathematica computer algebra system. This basic implementation sorts, performs row elimination of a matrix with entries in a polynomial ring with integer coefficients, such as $\mathbb{Z}[q,z]$, while not introducing rational functions, and it terminates when a row echelon matrix (a triangular system of equations) is reached. This function will be made a part of the impending next version release of \texttt{qFunctions} package. As a side note, we implemented a naive parallelization of this elimination but we have not seen any benefits of splitting calculations yet. 

We should also acknowledge that there are at least two crucial optimizations waiting to be implemented to aid proos of families in cylindric partitions scheme and other similar schemes. First task that should be done is to keep track of nullified relations and to remove the contributions of the nullspace in later calculations. To put it in concrete terms, at the moment we do not know if $N=6$ is the minimal number to prove Theorems~\ref{thm:n6m11} and/or \ref{thm:n6m13}. We know that it is a sufficient number. By removing any and all nullified relations we would only see a minimal representation (dependent on the choice of $N$) of these recurrences as elements in the ideals $I_m$, and that can give us an idea of what the optimal bound for $N$ is supposed to be in general. The second pending addition is dynamic extension of the matrix to be reduced. At the moment, we fix an $N$ experimentally hoping that it is enough to show that the relations of interest are in the nullspace of this matrix. This is in the same spirit of making a fixed ansatz. Row reduction as a preprocessing step helps for the repeated calculations. Having an echelon system boosts the speed of later calculations immensely. If the chosen $N$ is not enough, then we need to pick a larger $N$ and start all over. This requires performing row reduction of the matrix for $N$ once more as a subproblem. This should be changed by extending the already triangularized matrix for $N$ to $N+1$ and doing the row reduction again for only the added relations. The incrementality of the matrix would also carry us to the minimal necessary $N$ for any given $m$ (assuming that Conjecture~\ref{conj:FinIdeal} is correct) naturally.


\begin{thebibliography}{99}

\bibitem{qFunctions}
J.~Ablinger and A.~K. Uncu.
\newblock $\mathtt{qFunctions}$ - a {M}athematica package for $q$-series and
  partition theory applications.
\newblock Submitted. arXiv:1910.12410, 2019.

\bibitem{AAB87}
A.~Agrawal, G.~E. Andrews, and D.~Bressoud.
\newblock The {B}ailey lattice.
\newblock {\em J. Indian Math. Soc.}, 51:57--73, 1987.

\bibitem{And74}
G.~E. Andrews.
\newblock An analytic generalization of the {R}ogers-{R}amanujan identities for
  odd moduli.
\newblock {\em Proc. Nat. Acad. Sci. USA}, 71:4082--4085, 1974.


\bibitem{And86}
G. E. Andrews. 
\newblock \textit{q-series}: their development and application in analysis, number theory, combina-
torics, physics, and computer algebra. 
\newblock Vol. 66. CBMS Regional Conference Series in Math-
ematics. Published for the Conference Board of the Mathematical Sciences, Washington,
DC; by the American Mathematical Society, Providence, RI, 1986, pp. xii+130.

\bibitem{And84b}
G.~E. Andrews.
\newblock {\em The Theory of Partitions}.
\newblock Cambridge University Press, 1984.

\bibitem{And89}
G.~E. Andrews.
\newblock On the proofs of the {R}ogers-{R}amanujan identities.
\newblock In {\em $q$-Series and Partitions}, pages 1--14. Springer-Verlag, New
  York, 1989.

\bibitem{ASW99}
G.~E. Andrews, A.~Schilling, and S.~O. Warnaar.
\newblock An {$A_2$} {B}ailey lemma and {R}ogers-{R}amanujan-type identities.
\newblock {\em J. Amer. Math. Soc.}, 12(3):677--702, 1999.

\bibitem{AD11}
C. Armond and O. T. Dasbach. 
\newblock {R}ogers-{R}amanujan type identities and the head and tail of the colored {J}ones polynomial.
\newblock arXiv:1106.3948 [math.GT].

\bibitem{Bai49}
W.~N. Bailey.
\newblock Identities of the {R}ogers-{R}amanujan type.
\newblock {\em Proc. London Math. Soc.}, 50(2):1--10, 1949.

\bibitem{Bax81}
R.~J. Baxter.
\newblock {R}ogers-{R}amanujan identities in the hard hexagon model.
\newblock {\em J. Stat. Phys.}, 26:427--452, 1981.

\bibitem{Bor07}
A.~Borodin.
\newblock Periodic {S}chur process and cylindric partitions.
\newblock {\em Duke Math. J.}, 140(3):391--468, 2007.


\bibitem{BU} W. Bridges, and A.~K. Uncu.
\newblock {W}eighted cylindric partitions. 
\newblock {\em J. Algebraic Combin.}, 56 (2022), no. 4, 1309–-1337.

\bibitem{Bre79}
D.~M. Bressoud.
\newblock A generalization of the {R}ogers-{R}amanujan identities for all
  moduli.
\newblock {\em J. Comb. Th. A}, 27:64--68, 1979.

\bibitem{Bre83}
D.~M. Bressoud.
\newblock An easy proof of the {R}ogers-{R}amanujan identities.
\newblock {\em J. Number Th.}, 16:335--241, 1983.


\bibitem{BMS13}
C.~Bruschek, H.~Mourtada, and J.~Schepers.
\newblock Arc spaces and {R}ogers-{R}amanujan identities.
\newblock {\em Ramanujan J.}, 30:9--38, 2013.

\bibitem{Chern}
S.~Chern
\newblock {L}inked partition ideals, directed graphs and $q$-multi-summations.
\newblock {\em Electron. J, Combin.} 27(3): Paper No. 3.33, 29 pp.

\bibitem{Cor17}
S.~Corteel.
\newblock {R}ogers-{R}amanujan identities and the
  {R}obinson-{S}chensted-{K}nuth correspondence.
\newblock {\em Proc. Amer. Math. Soc.}, 145(5):2011--2022, 2017.

\bibitem{CW19}
S.~Corteel and T.~A. Welsh.
\newblock The {$A_2$} {R}ogers--{R}amanujan identities revisited.
\newblock {\em Annals of Combinatorics}, 23(3):683--694, 2019.

\bibitem{CDU}
S.~Corteel, J. Dousse and A.~K. Uncu.
\newblock {C}ylindric partitions and some new {A}2 {R}ogers–{R}amanujan identities. \newblock {\em Proc. Amer.Math. Soc.}, 150(2):481–-497, 2021.

\bibitem{FFW08}
B.~Feigin, O.~Foda, and T.~A. Welsh.
\newblock {A}ndrews--{G}ordon type identities from combinations of {V}irasoro
  characters.
\newblock {\em Ramanujan J.}, 17(1):33--52, 2008.

\bibitem{FW16}
O.~Foda and T.~A. Welsh.
\newblock Cylindric partitions, {$\mathcal{W}_r$} characters and the
  {A}ndrews-{G}ordon-{B}ressoud identities.
\newblock {\em J. Phys. A}, 49(16):164004, 37, 2016.

\bibitem{GM81}
A.~M. Garsia and S.~C. Milne.
\newblock A {R}ogers-{R}amanujan bijection.
\newblock {\em J. Combin. Theory Ser. A}, 31:289--339, 1981.

\bibitem{GK97}
I.~M. Gessel and C.~Krattenthaler.
\newblock Cylindric partitions.
\newblock {\em Trans. Amer. Math. Soc.}, 349(2):429--479, 1997.

\bibitem{Gor61}
B.~Gordon.
\newblock A combinatorial generalisation of the {R}ogers-{R}amanujan
  identities.
\newblock {\em Amer. J. Math.}, 83:393--399, 1961.

\bibitem{GOW16}
M.~J. Griffin, K.~Ono, and S.~O. Warnaar.
\newblock {A} framework of {R}ogers--{R}amanujan identities and their arithmetic
  properties.
\newblock {\em Duke Math. J.}, 8:1475--1527, 2016.

\bibitem{Kpriv}
S. Kanade.
\newblock {P}rivate communications.

\bibitem{KGG}
S. Kanade.
\newblock {S}tructure of certain level 2 standard modules for $A_5^{(2)}$ and {G}\"ollnitz--{G}ordon identities.
\newblock {\em Ramanujan J.}, 45(3):873--893, 2018.

\bibitem{KanF}
S. Kanade.
\newblock {O}n the {A}2 {A}ndrews–-{S}chilling–-{W}arnaar identities.
\newblock preprint.


\bibitem{KR}
S. Kanade, and M.~C. Russell.
\newblock {C}ompleting the {A}2 {A}ndrews--{S}chilling--{W}arnaar identities.
\newblock arXiv:2203.05690 [math.CO].

\bibitem{HolonomicFunctions}
C.~Koutschan.
\newblock {\em Advanced applications of the holonomic systems approach}.
\newblock PhD thesis, RISC, Johannes Kepler University, Linz, 2009.

\bibitem{LW84}
J.~Lepowsky and R.~L. Wilson.
\newblock The structure of standard modules, {I}: Universal algebras and the
  {R}ogers-{R}amanujan identities.
\newblock {\em Invent. Math.}, 77:199--290, 1984.

\bibitem{LW85}
J.~Lepowsky and R.~L. Wilson.
\newblock The structure of standard modules, {II}: The case ${A}_1^{(1)}$,
  principal gradation.
\newblock {\em Invent. Math.}, 79:417--442, 1985.

\bibitem{Mac16}
P.~A. MacMahon.
\newblock {\em Combinatory Analysis}, volume~2.
\newblock Cambridge University Press, New York, NY, USA, 1916.

\bibitem{ML92}
S.~C. Milne and G.~M. Lilly.
\newblock The {$A_{\ell}$} and {$C_{\ell}$} {B}ailey transform and lemma.
\newblock {\em Bull. Amer. Math.Soc.}, 26:258â€“263, 1992.

\bibitem{ML95}
S.~C. Milne and G.~M. Lilly.
\newblock Consequences of the {$A_{\ell}$} and {$C_{\ell}$} {B}ailey transform
  and lemma.
\newblock {\em Discrete Math.}, 139:319--346, 1995.

\bibitem{tinyurl}
A.K.~Uncu
\newblock $<$\url{https://drive.google.com/drive/folders/1qRLIfX8JVIzxkKQCCaYfg4i84X\_l\_-fo}$>$ {L}ast accessed January 3, 2023.

\bibitem{Pas20}
A.~Pascadi.
\newblock Several new product identities in relation to two-variable
  {R}ogers--{R}amanujan type sums and mock theta functions.
\newblock arXiv:2009.05878, 2020.

\bibitem{RR19}
L.~J. Rogers and S.~Ramanujan.
\newblock Proof of certain identities in combinatory analysis.
\newblock {\em Cambr. Phil. Soc. Proc.}, 19:211--216, 1919.

\bibitem{Sch17}
I.~Schur.
\newblock Ein {B}eitrag zur {A}dditiven {Z}ahlentheorie und zur {T}heorie der
  {K}ettenbr\"uche.
\newblock {\em S.-B. Preuss. Akad. Wiss. Phys. Math. Klasse}, pages 302--321,
  1917.

\bibitem{Sil17}
A.~V. Sills.
\newblock {\em An invitation to the Rogers-Ramanujan identities}.
\newblock CRC Press, 2017.

\bibitem{Tsu22}
S. Tsuchioka. 
\newblock {A}n example of {A}2 {R}ogers-{R}amanujan bipartition identities of level 3.
\newblock arXiv:2205.04811 [math.RT].

\bibitem{War21}
S.~O. Warnaar. \newblock {T}he {A}2 {A}ndrews-{G}ordon identities and cylindric partitions. 
\newblock arXiv:2111.07550 [math.CO].


\bibitem{WarPriv}
S.~O. Warnaar.
\newblock {P}rivate communications.


\end{thebibliography}
\end{document}